\documentclass[12pt]{amsart}
\usepackage{amssymb}
\usepackage{amsmath}
\usepackage{bbm}
\usepackage{dsfont}
\usepackage{epsfig}
\usepackage{graphicx}
\usepackage{mathptmx}
\usepackage{multicol}
\usepackage{multirow}
\usepackage{subfig}
\usepackage{verbatim}
\usepackage{pgf,tikz}
\usepackage{tkz-fct}
\usepackage{pgfplots}
\usetikzlibrary{arrows,trees}
\definecolor{darkergreen}{RGB}{0,153,0}
\definecolor{darkerblue}{RGB}{0,0,204}

\newcounter{teorema}
\newenvironment{teorema}{\goodbreak
   \refstepcounter{teorema}
   \medskip\noindent{\bf Theorem~\theteorema .}\hskip 5pt}{}
\newenvironment{lemma}{\goodbreak
   \refstepcounter{teorema}
   \medskip\noindent{\bf Lemma~\theteorema .}\hskip 5pt}{}

\newenvironment{proposition}{\goodbreak
   \refstepcounter{teorema}
   \medskip\noindent{\bf Proposition~\theteorema .}\hskip 5pt}{}

\def\E{\mathbb{E}}
\def\N{\mathbb{N}}
\def\R{\mathbb{R}}

\def\Z{\mathbb{Z}}

\def\qad{\hskip 5pt}
\def\mod#1{\,({\rm mod\ }#1) }
\def\proof{\goodbreak\bigskip\noindent {\sc Proof}. \/}
\begin{document}

\title{Nonlinear rotations on a lattice}
\author{Fairuz Alwani}
\author{Franco Vivaldi}
\address{School of Mathematical Sciences, Queen Mary,
University of London,
London E1 4NS, UK}

\begin{abstract} We consider a prototypical two-parameter family of invertible maps 
of $\mathbb{Z}^2$, representing rotations with decreasing rotation number.
These maps describe the dynamics inside the island chains of
a piecewise affine discrete twist map of the torus, in the limit of fine discretisation.
We prove that there is a set of full density of points which, depending of the parameter 
values, are either periodic or escape to infinity. 
The proof is based on the analysis of an interval-exchange map over the integers,
with infinitely many intervals.
\end{abstract}
\maketitle
\section{Introduction}\label{section:Introduction}

Regular motions in two-dimensional symplectic maps are rotations on smooth invariant curves. 
If the space is discrete (a lattice, typically), then these curves do not exist, intriguing 
new phenomena appear, and the stability problem ---the central problem of Hamiltonian 
mechanics--- must be reconsidered from scratch.

Discrete-space versions of symplectic maps first appeared in the study of numerical 
orbits \cite{Rannou,Kaneko,Scovel,EarnTremaine,Vivaldi:94,NucinkisEtAl}, 
to mimic quantum effects in classical systems \cite{ChirikovEtAl}, and to 
improve the efficiency of delicate computations \cite{Karney}. 
The arithmetical characterisation of chaotic orbits provided a new direction of research
\cite{PercivalVivaldi,Keating,DegliEspostiIsola,NeumaerkerRobertsVivaldi},
and so did the study of the dynamics of round-off errors
\cite{LowensteinHatjispyrosVivaldi,LowensteinVivaldi:98,LowensteinVivaldi:00,
BosioVivaldi,VivaldiVladimirov,KouptsovLowensteinVivaldi,ReeveBlackVivaldi:14}.
Discrete symplectic maps occur in the study of outer billiards of polygons
\cite{Schwartz}, and in shift-radix systems in arithmetic
\cite{AkiyamaBrunottePethoThuswaldner,AkiyamaBrunottePethoSteiner};

In spite of a protracted research effort, our knowledge of these systems remains
fragmented; in particular, the stability problem has proved stubbornly difficult.
Rigorous results are rare, and the many and varied mechanisms responsible for 
(in)stability do not yet fit into a coherent picture, let alone a mathematical
theory.

\begin{figure}[t]
\hspace*{30pt}
\epsfig{file=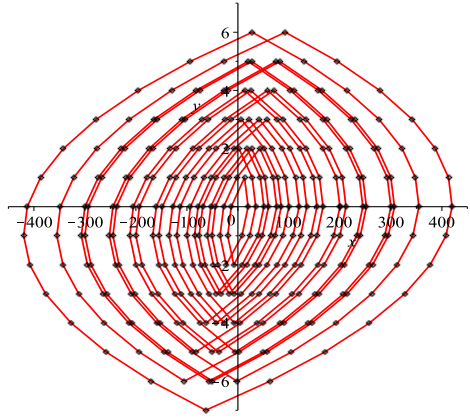,width=12cm,height=14cm}
\vspace*{-210pt}
\caption{\label{fig:LatticeMap}\rm\small
An orbit of the map (\ref{eq:LatticeMap}), with $\alpha=19$ and $\beta=7$.
In spite of large fluctuations in amplitude, the orbit closes up after $\alpha$ 
revolutions around the origin. Nearby orbits are intertwined, and hence their
boundedness cannot be inferred from topological considerations.
}
\end{figure}

If an area-preserving map preserves a lattice, then the existence of 
bounding invariant sets for the map guarantees the boundedness of lattice orbits.
Examples include the invariant polygons of the saw-tooth map \cite{Devaney}, and 
the invariant necklaces in outer billiards of (quasi)-rational polygons \cite{VivaldiShaidenko}.
But bounding invariant sets in an embedding space of a lattice are rarely available,
and a different approach is needed. In the case of rational rotations on lattices
(which necessarily involve some rounding procedure),
all available proofs of stability rely on renormalization, which provides knowledge of
long-time asymptotics \cite{LowensteinHatjispyrosVivaldi,
KouptsovLowensteinVivaldi,AkiyamaBrunottePethoThuswaldner,AkiyamaBrunottePethoSteiner}.
Renormalization was also key to the proof of the existence of escape orbits 
of outer billiards of kites over quadratic fields \cite{Schwartz:07}.
However, renormalizability too is seldom available (for rational rotations on a
lattice it occurs only for finitely many quadratic irrational parameter values).
Thus no proof of stability is known for invertible irrational rotations on lattices,
even though the orbits are believed to be periodic.
Here the round-off perturbation generates diffusive transport, yet all orbits seemingly 
return to their initial point via a mechanism that is probabilistic at heart
\cite{BosioVivaldi,VivaldiVladimirov}.
The observed stability of the rotational orbits of certain linked strip maps on lattices
is even more elusive \cite{ReeveBlackVivaldi:14}.

In this paper we illustrate a novel mechanism for the (in)stability of rotational orbits,
as it appears in the following two-parameter family of invertible nonlinear maps $\mathcal{F}$ of 
the two-dimensional lattice $\mathbb{Z}^2$ 
\begin{equation}\label{eq:LatticeMap}
\begin{array}{rcl}
y_{t+1}&=&y_t-\mbox{sign}(x_t)\\
x_{t+1}&=&x_t+\alpha y_{t+1}+\beta
\end{array}
\qquad  
\mbox{sign}(x)=
\begin{cases}
1 & \mbox{if}\,\,\, x\geqslant 0\\
-1& \mbox{if}\,\,\, x<0,
\end{cases}
\end{equation}
where $\alpha$ and $\beta$ are integers, and $0\leqslant \beta <\alpha$.
(There is no loss of generality in choosing a perturbation of unit magnitude, 
since the perturbation amplitude can be absorbed by the other parameters.)
As is often the case in piecewise affine dynamical systems, the plain form of 
(\ref{eq:LatticeMap}) hides a non-trivial dynamics (see figure \ref{fig:LatticeMap}); 
concatenated parabolic arcs result in surrogate rotations with decreasing rotation number.

The map $\mathcal{F}$ originates from the following perturbed twist map on a discrete torus \cite{ZhangVivaldi}:
\begin{equation}\label{eq:TorusMap}
\begin{array}{rcll}
y_{t+1}&\equiv&y_t+f(x_t) &\mod{N}\\
x_{t+1}&\equiv&x_t+y_{t+1}&\mod{N}
\end{array}
\hskip 40pt 
f(q)=\begin{cases}1 & 0\leqslant q< \lfloor N/2\rfloor\\
-1 & \mbox{otherwise.}
\end{cases}
\end{equation}
Here $N$ is a large integer ---the discretisation parameter---
while the perturbation function $f$ provides a minimalist form of nonlinearity.
This map is a `pseudo-elliptic' variant of the so-called triangle map,
which is `pseudo-hyperbolic'. (The puzzling ergodic properties of the latter have 
so far escaped a rigorous analysis \cite{CasatiProsen,HorvatEtAl,NeumaerkerRobertsVivaldi}.)

\begin{figure}[t]
\hspace*{20pt}
\epsfig{file=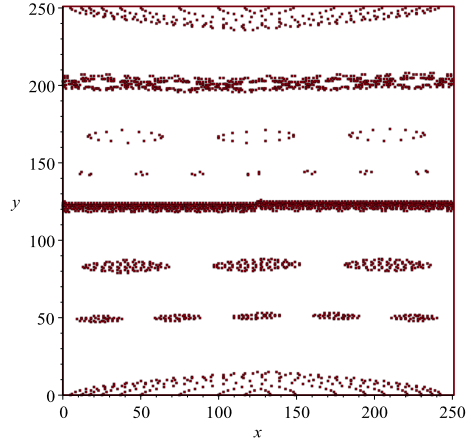,width=12cm,height=15cm}
\vspace*{-220pt}
\caption{\label{fig:TorusMap}\rm\small
Some orbits of the perturbed twist mapping (\ref{eq:TorusMap}), for $N=251$.
Asymptotic ($N\to\infty$) dynamics inside island chains are described by the 
map $\mathcal{F}$ given in (\ref{eq:LatticeMap}). 
}
\end{figure}

In figure \ref{fig:TorusMap} we display some orbits of the map (\ref{eq:TorusMap}),
which bears resemblance to the divided phase space of an area-preserving map.
We observe island chains of odd order, those of even order are missing, and
there is no hierarchy of islands about islands.
Plainly, standard Hamiltonian perturbation theory does not apply, so
what does determine the stability of these elliptic orbits?

It can be shown that, for sufficiently large $N$, the map (\ref{eq:LatticeMap}) is the 
first-return map to an island of (\ref{eq:TorusMap}) for all points sufficiently close 
to the island's centre. If the island has rotation number $m/n$, then $\alpha=n$, and 
$$
\beta=\beta(m,n,N)=\sum_{t=1}^{n-1}(-1)^{\lfloor 2mt/n\rfloor}-mN \mod{n}.
$$
(See \cite{ZhangVivaldi} for details.)

In this paper we solve the stability problem of $\mathcal{F}$ for a set initial conditions 
having full density. We prove that, depending on the parameters values, the orbits are either
periodic, or escape to infinity in both time directions, as accelerator modes.
In the final analysis, (in)stability will result from ergodicity in an associated 
modular arithmetic system.

\begin{figure}[t]
\hspace*{30pt}
\epsfig{file=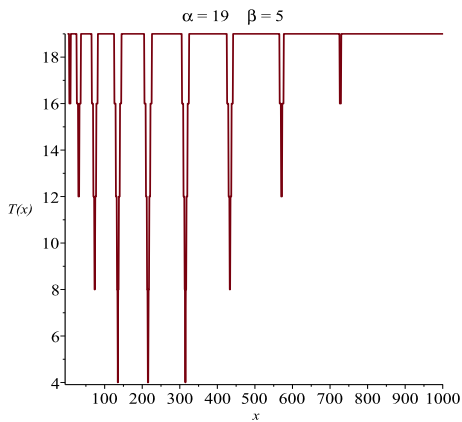,width=12cm,height=14cm}
\vspace*{-210pt}
\caption{\label{fig:Periods}\rm\small
Period $T(x)$ of the orbit though $x$ for the interval-exchange map $\mathrm{F}$ associated with (\ref{eq:LatticeMap}),
with $\alpha=19$ and $\beta=5$. (The vertical segments in the graph of $T$ are merely a guide to the eye.)
The behaviour near the origin is complicated, but for sufficiently large initial points 
($x\geqslant 730$) the period stabilises at $\alpha$.
The depth and width of this comb-like structure depends sensitively on the
arithmetical properties of the parameters.
}
\end{figure}

We will show that the first-return map $\mathrm{F}$ to the ray 
$\{(x,0)\in\Z^2\,:\, x\geqslant 0\}$ is an interval-exchange transformation
over infinitely many intervals. Near the origin, the dynamics is rather 
intricate (see figure \ref{fig:Periods}), but at large amplitudes, 
the map $\mathrm{F}$ admits a weak form of translational invariance. 
The large-amplitude dynamics is captured by the following conjecture (cf.~\cite{ZhangVivaldi})

\bigskip \noindent {\bf Conjecture.} 
{\sl Let $\overline{\alpha}=\alpha/\gcd(\alpha,2\beta)$.
If $\overline{\alpha}$ is odd, then all orbits of $\mathcal{F}$ are periodic,
and for all but finitely many initial conditions, their period under the 
first-return map $\mathrm{F}$ is equal to $\overline{\alpha}$. If $\overline{\alpha}$ is
even, then all orbits escape to infinity.
}
\medskip

This conjecture is consistent with the absence of island chains or even order,
observed experimentally for the map (\ref{eq:TorusMap}).
The main result of this paper is the following theorem, which establishes 
a probabilistic version of the above conjecture.

\begin{teorema}\label{thm:Main}
{\sl 
If $\overline{\alpha}$ (as defined above) is odd, then the periodic points of
$\mathcal{F}$ have full density, and their period under the first-return map 
$\mathrm{F}$ is equal to $\overline{\alpha}$.
If $\overline\alpha$ is even, then the set of escape orbits has full density.
}
\end{teorema}
\medskip

The first-return map $\mathrm{F}$ will be constructed in section \ref{section:FirstReturnMap},
where we derive several formulae to be used throughout the paper.
In section \ref{section:Parameters} we show that there is no loss of generality in 
restricting the parameters to the range $\alpha\geqslant 2\beta$ with $\alpha$ and $\beta$ co-prime 
(propositions \ref{thm:ParameterSpecialCases}--\ref{thm:ParameterConjugacy}). 
In section \ref{section:IET} we show that $\mathrm{F}$ is an interval-exchange
transformation over infinitely many intervals; we compute the IET's metric data,
and establish that the combinatorial data are (essentially) parameter-independent 
(proposition \ref{prop:sigma}).

In section \ref{section:SymbolicDynamics} we consider the natural symbolic dynamics
of the IET, together with two coarser codes, to factor out
translations in the code, and to anchor the code to the minimum point of an orbit.
Asymptotically, the cylinder sets of the symbolic dynamics have a regular 
structure ---they are arranged into arithmetic progressions.

In section \ref{section:ReducedSystem} we derive an auxiliary interval-exchange 
map $\mathrm{F}'$ over $\mathbb{Z}$ ---the \textit{reduced system}--- which 
encodes the asymptotic behaviour of the original IET.
The idea is to take the large-amplitude limit of $\mathrm{F}$, scale it 
in such a way as to obtain a spatially periodic integer map, and then 
extend the latter periodically to $\Z$. The periodic cells of the reduced system
are the \textit{blocks}, the union of two adjacent intervals of the IET.
We prove that our conjecture holds for the reduced system (theorem 
\ref{thm:ReducedSystem}).

In section \ref{section:RegularPoints} we consider the \textit{regular points} of the 
Poincar\'e map $\mathrm{F}$, namely the points whose symbolic words of length $\alpha$ 
also belong to the language of the reduced map $\mathrm{F}'$.
We then prove that almost all points are regular (theorem \ref{thm:FullDensity}),
which will allow us to use the symbolic dynamics of the reduced system for the
original system.

Theorem \ref{thm:Main} is proven is sections \ref{section:PeriodicOrbits} 
and \ref{section:EscapeOrbits}.
To establish the periodicity of all regular points of $\mathrm{F}$,
we must determine the value of a certain invariant of the reduced system. 
This invariant behaves like a variance, and the key lemma 
\ref{lemma:Key} establishes its value by considering the evolution of uniform
measures supported on blocks.
A similar technique is used in section \ref{section:EscapeOrbits}, 
to show that, if $\overline\alpha$ is even, then almost all orbits escape.
In this case however, the aforementioned invariant is replaced by a non-constant
function of the coordinates, whose regular variation is determined using the
Sturmian property of rotational codes.
The computations of this section are considerably more laborious than for the periodic case.

The map (\ref{eq:LatticeMap}) admits natural generalisations to higher-dimensional
lattices. For instance, one could choose the parameters $\alpha$ and $\beta$ from some ring 
$\mathbb{Z}[\omega]$ of real algebraic integers, to obtain a dynamical system over 
$\mathbb{Z}[\omega]^2$ (or, more generally, over the Cartesian product of
two $\mathbb{Z}[\omega]$-modules). These are four-dimensional lattices, and there
is no reason to expect theorem \ref{thm:Main} to extend to such systems. 
In numerical experiments over quadratic fields, we have observed recurrence and a 
weak form of instability replacing periodicity.

\section{First-return map}\label{section:FirstReturnMap}

In this section we construct the first-return map $\mathrm{F}$ to the ray 
$\mathbb{Z}_+=\{(x,0)\,:\,x\geqslant 0\}$, which is crossed repeatedly by 
every orbit of $\mathcal{F}$. 
Let $\mathbb{Z}_-=\{(x,0)\,:\,x< 0\}$. To construct $\mathrm{F}$, we consider the first 
transit maps $\mathrm{F}_\pm$ from $\Z_\pm$ to $\Z$:
$$
\mathrm{F}_+:\Z_+\to\Z,
\hskip 40pt
\mathrm{F}_-:\Z_-\to\Z.
$$
The idea is to define $\mathrm{F}=\mathrm{F}_-\circ \mathrm{F}_+$. This is legitimate only if
$\mathrm{F}_+$ and $\mathrm{F}_-$ map $\Z_+$ to $\Z_-$, and vice-versa. 
As we shall see, this is not always the case.

We begin by solving (\ref{eq:LatticeMap}) over each domain where $\mbox{sign}(x)$ remains constant.
Specifically, let $x_0$ and $t$ be such that $\mbox{sign}(x_k)=\mbox{sign}(x_0)$ for $k=0,\ldots,t$. 
We compute:
\begin{equation}\label{eq:Solution}
y_t=y_0-st\qquad x_t=x_0-\frac{t(t+1)}{2}\alpha s+t(\alpha y_0+\beta)
\qquad s=\mbox{sign}(x_0).
\end{equation}
Let now $u_s$ ($s$ as above) be the smallest positive integer $t$ such that $\mbox{sign}(x_t)\not=s$.
To construct $\mathrm{F}_\pm$, we specialise formula (\ref{eq:Solution}) to the initial conditions 
$(x_0,y_0)=(x,0)$, and then match two solutions (\ref{eq:Solution}) near $x=0$, to obtain
\begin{equation}\label{eq:Fpm}
\begin{array}{rcll}
\mathrm{F}_+(x)&=&x+\tau_+(x)&\qquad x\geqslant 0\\
\mathrm{F}_-(x)&=&x+\tau_-(x)&\qquad x<0
\end{array}
\end{equation}
where
\begin{eqnarray}
\tau_+(x)&=&2\beta u_+(x)-\alpha u_+(x)^2\nonumber \\
u_+(x)&=&\lfloor U_+(x)+1\rfloor\nonumber\\
U_+(x)&=&\displaystyle \frac{1}{2\alpha}\left(2\beta-\alpha
   +\sqrt{(2\beta-\alpha)^2+8\alpha x}\right) \nonumber \\
\nonumber \\
\tau_-(x)&=&2\beta u_-(x)+\alpha u_-(x)^2\nonumber \\
u_-(x)&=&\lceil U_-(x)\rceil\nonumber\\
U_-(x)&=&\displaystyle \frac{1}{2\alpha}\left(-(2\beta+\alpha)
   +\sqrt{(2\beta+\alpha)^2-8\alpha x}\right). \label{eq:Um}
\end{eqnarray}
As functions over $\R$, the functions $\tau_\pm$ are singular, and right continuous at
each singularity, as easily verified.

As $x\to\infty$, we have  $\tau_+(x)\sim -2x$; 
likewise as $x\to -\infty$ we have $\tau_-(x)\sim 2x$.
Hence, for all sufficiently large $x$, we have 
$\mathrm{F}_+(x)<0$ and $\mathrm{F}_-(-x)>0$.

Next we investigate what happens for small $x$.
Let $(x_m)$ and $(y_n)$ be the sequences of singularities of 
$\mathrm{F}_+$ and $\mathrm{F}_-$, respectively. We compute
\begin{equation}\label{eq:xmyn}
\begin{array}{rcll}
x_m=\displaystyle\frac{m}{2}\left(\alpha(m+1)-2\beta\right)
   &\quad m=0,1,2,\ldots\\
\noalign{\vskip 3pt}
y_n=-\displaystyle\frac{n}{2}\left(\alpha(n+1)+2\beta\right)
   &\quad n=0,1,2,\ldots
\end{array}
\end{equation}
where the case $n=0$ is introduced for convenience.

We begin with $\mathrm{F}_+$. Letting
\begin{equation}
z_m=\mathrm{F}_+(x_m)=x_m+\tau_+(x_m)\nonumber
\end{equation}
we find
\begin{eqnarray}
z_m&=&\displaystyle\frac{m}{2}(\alpha (m+1)-2\beta)+(m+1)(2\beta-\alpha(m+1))\nonumber\\
  &=&\displaystyle\frac{m+2}{2}\bigl(2\beta-\alpha (m+1)\bigr).\label{eq:zm}
\end{eqnarray}
Note that $x_m,y_n$ and $z_m$ are integers, and that
\begin{equation}\label{eq:zxIdentity}
x_{m+1}-x_m=z_m-z_{m-1}=\alpha(m+1)-\beta\qquad m\geqslant 1.
\end{equation}

Now, since $\tau_+$ is right-continuous, for all $m\geqslant 0$ we have
\begin{equation}\label{eq:Fp}
y=\mathrm{F}_+(x)=x-x_m+z_m\hskip 40pt x_m\leqslant x <x_{m+1}
\end{equation}
and we find
\begin{equation}
a_m=\sup_{x_m\leqslant x <x_{m+1}}\mathrm{F}_+(x) =x_{m+1}-x_m+z_m
   =\frac{m+1}{2}(2\beta-\alpha m).\nonumber
\end{equation}
Thus $a_m<0$ for $m>1$, whereas $a_m>0$ for $m=0$, and, if 
$\alpha<2\beta$, also for $m=1$.

We repeat the analysis for $\mathrm{F}_-$. 
We define \begin{equation}\label{eq:wm}
w_n=\mathrm{F}_-(y_n)=\frac{n}{2}(\alpha(n-1)+2\beta)\qquad n\geqslant 1.
\end{equation}
Now, since $\tau_-$ is right-continuous, we have
\begin{equation}\label{eq:Fm}
w=\mathrm{F}_-(y)=y-y_n+w_n\hskip 40pt y_n\leqslant y <y_{n-1}
\end{equation}
so that
\begin{equation}
b_n=\min_{y_n\leqslant y <y_{n-1}}\mathrm{F}_-(y)=w_n.\nonumber
\end{equation}
Thus $b_n\geqslant 0$ for all $n$.

We now extend the domain of $\mathrm{F}_-$ to include all positive values 
of $x$ for which $U_-(x)$ is real ---see (\ref{eq:Um}).
We find that
\begin{equation}\label{eq:GoodRange}
-1<U_-(x)\leqslant 0
\qquad \mbox{if}\qquad 
0\leqslant x\leqslant x^*:=\frac{(2\beta+\alpha)^2}{8\alpha}
\end{equation}
so that in this $x$-range we have $u_-(x)=0$ and $\mathrm{F}_-(x)=x$. 
We verify that the image of $\mathrm{F}_+$ remains within this range:
$$
x^*-a_0=\frac{(\alpha-2\beta)^2}{8\alpha}\geqslant 0
\hskip 40pt
x^*-a_1=\frac{(3\alpha-2\beta)^2}{8\alpha}\geqslant 0.
$$
So we have
$$
\mathrm{F}_+(x)\geqslant 0\quad \Rightarrow \quad \mathrm{F}_-(\mathrm{F}_+(x))=\mathrm{F}_+(x),
$$
hence $\mathrm{F}=\mathrm{F}_-\circ \mathrm{F}_+$, and we have established the following result:

\begin{proposition} \label{thm:F}
{\sl The first-return map $\mathrm{F}$ to $\Z_+$ is of the form $x\mapsto x+\tau(x)$},
where
\begin{equation*}\label{eq:F}
\tau(x)=2\beta (u_-(y)+u_+(x))+\alpha (u_-^2(y)-u_+^2(x))\hskip 30pt y=\mathrm{F}_+(x).
\end{equation*}
\end{proposition}
\bigskip

Next we compute the sequences of singularities of $\mathrm{F}$. 
To this end, we must determine 
the sequence $(\mathrm{F}_+^{-1}(y_n))$ and then merge it with $(x_m)$.
To compute $\mathrm{F}_+^{-1}$ we solve (\ref{eq:Fp}) for $x$, and then use
(\ref{eq:xmyn}) and (\ref{eq:zm}) to obtain
\begin{equation}\label{eq:F+inverse}
\mathrm{F}_+^{-1}(y)=y+(m(y)+1)(\alpha (m(y)+1)-2\beta)
\hskip 30pt
z_m\leqslant y < z_{m-1}.
\end{equation}
Here $m$ is the smallest integer such that $z_m\leqslant y$. 
We find
$$
m(y)=
\left\lceil
\frac{2\beta-3\alpha+\sqrt{(2\beta+\alpha)^2-8\alpha y}}{2\alpha}
\right\rceil
$$
and one verifies that $ m(y)+1=\lceil U_+(\beta-y)\rceil.$

To order the singularities of $\mathrm{F}$ we must establish a relationship 
between the indices $m$ and $n$, namely find all solutions $n=n(m)$ 
of the inequalities
$$
z_m\leqslant y_n<z_{m-1}.
$$
In what follows we exclude the special cases $\beta=0$ and $\alpha=2\beta$ 
which are dealt with in proposition \ref{thm:ParameterSpecialCases} in the next section.
Let $n=m+k$. The lower and upper bounds give, respectively
\begin{eqnarray}
i)&&((k-1)\alpha+2\beta)(k+2m+2)\leqslant 0\nonumber\\
ii)&&(k\alpha+2\beta)(k+2m+1) > 0. \nonumber
\end{eqnarray}
We obtain
\begin{eqnarray}
i)&&-2(m+1)\leqslant k\leqslant 1-\frac{2\beta}{\alpha}\nonumber \\
ii)&&k<-(2m+1)\quad\mbox{or}\quad k>-\frac{2\beta}{\alpha}. \nonumber
\end{eqnarray}
Since $m+k\geqslant 0$, the relevant bound in $ii)$ is the rightmost one, and
we find
\begin{equation}\label{eq:k}
k=
\begin{cases}
 0&\alpha >2\beta\\
-1&\alpha < 2\beta.
\end{cases}
\end{equation}
Accordingly, we let
$$
x'_m=
\begin{cases}
\mathrm{F}_+^{-1}(y_m)    &\alpha>2\beta\\ 
\mathrm{F}_+^{-1}(y_{m-1})&\alpha< 2\beta 
\end{cases}
$$
where in both cases we use the same branch of $\mathrm{F}_+^{-1}$,
specified in (\ref{eq:F+inverse}).
We find:
\begin{equation}\label{eq:xpm}
x_m'=
\begin{cases}
\frac{m}{2}\left(\alpha m+3(\alpha-2\beta)\right)
    +\alpha-2\beta   &\alpha>2\beta \\
\frac{m}{2}\left(\alpha m+5\alpha-6\beta\right)
    +\alpha -\beta&\alpha< 2\beta.
\end{cases}
\end{equation}

Let $(\delta_m)$, $m\geqslant 0$ be the sequence of singularities 
of $\mathrm{F}$, in ascending order. 
From (\ref{eq:k}) we have $\delta_0=0$ and
\begin{equation}\label{eq:delta}
\delta_{2m}=\begin{cases} 
 x_m &\alpha>2\beta \\
 x_m'&\alpha<2\beta 
\end{cases}
\qquad
\delta_{2m-1}=\begin{cases} 
 x_{m-1}' &\alpha>2\beta \\
 x_{m}&\alpha<2\beta 
\end{cases}
\qquad m=1,2,\ldots.
\end{equation}
This leads to the sequences of singularities
$$
\begin{array}{lll}
(x_0,x_0',x_1,x_1',x_2,x_2',\ldots) && \alpha>2\beta\\
\noalign{\vskip 1pt}
(x_0,x_1,x_1',x_2,x_2',x_3,\ldots)  && \alpha<2\beta.
\end{array}
$$

\section{Parameters}\label{section:Parameters}

In this section we show that there is no loss of generality in restricting 
the parameters of the map $\mathcal{F}$ to the 
range $\alpha>2\beta$ with $\alpha$ and $\beta$ co-prime.
This is the content of the following three propositions.
To make the parameter-dependence of $\mathcal{F}$ and $\mathrm{F}$ explicit, we 
shall use the notation ${\mathcal F}_{\alpha,\beta}$ and $\mathrm{F}_{\alpha,\beta}$. 

First, we dispose of the special parameter values $\beta=0$ and $2\beta=\alpha$,
at which the singularities of $\mathrm{F}$ cancel out and the dynamics is trivial.

\begin{proposition} \label{thm:ParameterSpecialCases}
{\sl
If $\beta=0$ or $2\beta=\alpha$, then $\mathrm{F}$ is the identity.
}
\end{proposition}

\medskip

\proof
Let $\beta=0$. From (\ref{eq:xmyn}), (\ref{eq:zm}), and (\ref{eq:wm}) 
we verify that $z_m=y_{m+1}$ and that $w_{m+1}=x_m$. We find
$$
x_m=w_{m+1}=\mathrm{F}_-(y_{m+1})=\mathrm{F}_-(z_{m})=\mathrm{F}_-(\mathrm{F}_+(x_{m}))=\mathrm{F}(x_{m}).
$$
Likewise, if $2\beta=\alpha$, then $z_m=y_m$ and $w_m=x_m$, and we have
$$
x_m=w_{m}=\mathrm{F}_-(y_{m})=\mathrm{F}_-(z_{m})=\mathrm{F}_-(\mathrm{F}_+(x_{m}))=\mathrm{F}(x_{m}).
$$
In both cases the sequences $(x_m')$ and $(x_m)$ map into one another,
and $x_m$ is a fixed point of $\mathrm{F}$ for all $m$.

Our claim now follows from the fact that the function $\tau$ is
piecewise-constant and right-continuous at all its singularities.\quad $\Box$.

\bigskip
Next we reduce the size of parameter space by establishing a symmetry. 

\begin{proposition} \label{thm:ParameterSymmetry}
{\sl
For all $\alpha,\beta$ we have $\mathrm{F}_{\alpha,\alpha-\beta}=\mathrm{F}^{-1}_{\alpha,\beta}$.
}
\end{proposition}

\medskip
\proof  
First we show that the singularities of the two maps coincide. 
The singularities of $\mathrm{F}_{\alpha,\beta}^{-1}$ are $\mathrm{F}_{\alpha,\beta}(x'_m)$
and $\mathrm{F}_{\alpha,\beta}(x_m)$. 
We shall use equations (\ref{eq:xmyn}--\ref{eq:Fm}).
For all $\alpha,\beta$ we have
\begin{eqnarray}
\mathrm{F}_{\alpha,\beta}(x_m'(\alpha,\beta))&=&\mathrm{F}_-(y_m)=w_m=\frac{m}{2}(\alpha(m-1)+2\beta)\nonumber\\
 &=&x_m(\alpha,\alpha-\beta). 
\label{eq:Singularities}
\end{eqnarray}
For $\alpha > 2\beta$ we have
\begin{eqnarray}
\mathrm{F}_{\alpha,\beta}(x_m(\alpha,\beta))&=&\mathrm{F}_-(z_m)=z_m-y_{m+1}+w_{m+1}\label{eq:SingularitiesII}\\
&=&
\frac{m}{2}\left(\alpha(m+1)+6\beta)\right) +4\beta\nonumber\\
&=&x_{m+1}'(\alpha,\alpha-\beta).\nonumber
\end{eqnarray}
For $\alpha < 2\beta$ we have
\begin{eqnarray}
\mathrm{F}_{\alpha,\beta}(x_m(\alpha,\beta))&=&\mathrm{F}_-(z_m)=z_m-y_m+w_m\label{eq:SingularitiesIII}\\
&=&
\frac{m}{2}\left(\alpha m-3(\alpha-2\beta)\right)
    -(\alpha-2\beta)\nonumber\\
&=&x_m'(\alpha,\alpha-\beta).\nonumber
\end{eqnarray}
So $\mathrm{F}_{\alpha,\beta}$ and $\mathrm{F}^{-1}_{\alpha,\alpha-\beta}$ have the same singularities.
This result, together with the analogous calculations with exchanged parameters,
show that the value of $\mathrm{F}_{\alpha,\beta}$ and $\mathrm{F}_{\alpha,\alpha-\beta}^{-1}$ 
at those singularities is the same. The right continuity of the functions $\mathrm{F}$ and $\mathrm{F}^{-1}$ 
establishes the result. \/ $\Box$

\medskip

Finally, we show that it suffices to consider the case $\gcd(\alpha,\beta)=1$.
Let $d$ be a positive integer, and let us
consider the map ${\mathcal L}_{d\alpha,d\beta}$, with $\gcd(\alpha,\beta)=1$.
Then, for any $r$ in the range $0\leqslant r<d$, the set
\begin{equation}\label{eq:Ldr}
\mathbb{L}_{d,r}=(r+d\Z)\times \Z
\end{equation}
is invariant under ${\mathcal L}_{d\alpha,d\beta}$ [since in this case
$x_{t+1}\equiv x_t \mod{d}$, from (\ref{eq:LatticeMap})].

\begin{proposition} \label{thm:ParameterConjugacy} {\sl
Let $d\in\N$.
Then, for any $r$ in the range $r\in\{0,\ldots,d-1\}$, the map ${\mathcal F}_{\alpha,\beta}$ 
is conjugate to the restriction of ${\mathcal F}_{d\alpha,d\beta}$ to $\mathbb{L}_{d,r}$.
}
\end{proposition}

\proof The map
$$
\psi_r:\Z^2\to (r+d\Z)\times \Z\hskip 40pt (x,y)\mapsto (r+dx,y)
$$
is clearly a bijection. We must show that
$$
\psi_r\circ{\mathcal F}_{\alpha,\beta}=
\left.{\mathcal F}_{d\alpha,d\beta}\right|_{{\mathbb L}_{d,r}} \circ \psi_r.
$$
We compute
\begin{eqnarray}
(\psi_r\circ{\mathcal F}_{\alpha,\beta})(x,y)&=&
\psi_r(x+\alpha y -\alpha \mbox{sign}(x)+\beta,y-\mbox{sign}(x))\nonumber\\
&=&(r+d(x+\alpha y -\alpha \mbox{sign}(x)+\beta),y-\mbox{sign}(x))\label{eq:Conjugacy}\\
&=&(r+dx+d\alpha y -d\alpha\mbox{sign}(x)+d\beta),y-\mbox{sign}(x)).\nonumber
\end{eqnarray}
Now, for any $x\in\Z$ we have $\mbox{sign}(x)=\mbox{sign}(dx+r)$.
This is clearly true if $x\geqslant 0$, since $r\geqslant 0$.
If $x<0$, then
$$
dx+r\leqslant -d+r\leqslant -1
$$
and hence $dx+r$ has the same sign as $x$.
Using this identity in (\ref{eq:Conjugacy}), we obtain:
\begin{eqnarray*}
(\psi_r\circ{\mathcal F}_{\alpha,\beta})(x,y)
&=& \left[(r+dx)+d\alpha y -d\alpha\mbox{sign}(r+dx)+d\beta,y-\mbox{sign}(r+dx)\right]\\
&=&
\left(\left.{\mathcal F}_{d\alpha,d\beta}\right|_{{\mathbb L}_{d,r}} \circ \psi_r\right)(x,y),
\end{eqnarray*}
as desired.
$\Box$

\section{The interval-exchange transformation}\label{section:IET}

In this section we characterise the first-return map $\mathrm{F}$ defined
in section \ref{section:FirstReturnMap} as an interval-exchange 
transformation, by computing its metric and combinatorial data.
There are only two distinct permutations of the intervals, corresponding 
to the two parameter ranges $\alpha>2\beta$ and $\alpha<2\beta$,
one permutation being the inverse of the other (proposition \ref{prop:sigma}).

We define the sequence of intervals
\begin{equation}\label{eq:Delta}
\Delta_m=[\delta_{m-1} ,\delta_{m})\quad m=1,2,\ldots
\end{equation}
where $\delta_m$ is defined in (\ref{eq:delta}). 
These intervals form a partition of $\Z_+$.
The restriction of $\mathrm{F}$ to each interval is a translation, and hence 
$\mathrm{F}$ ---being invertible--- is an interval-exchange transformation.

For $\alpha>2\beta$, and $m\geqslant 1$ the corresponding translations are given by
$$
\tau_{2m}= \mathrm{F}(x_{m-1}')-x_{m-1}',
\qquad
\tau_{2m-1}= \mathrm{F}(x_{m-1})-x_{m-1}
$$
while the interval lengths are
$$
|\Delta_{2m}|= x_{m}-x_{m-1}',
\hskip 30pt
|\Delta_{2m-1}|= x_{m-1}'-x_{m-1}.
$$
Using (\ref{eq:Singularities}) we obtain 
\begin{equation}\label{eq:Parameters}
\begin{array}{rcll}
\tau_{2m}&=& (2\beta-\alpha)(2m-1)\\
\tau_{2m-1}&=& 4\beta m
&\qquad\multirow{2}{100pt}{$\alpha>2\beta, \hskip 3pt m\geqslant 1$}\\
|\Delta_{2m}|&=& \beta(2m-1)\\
|\Delta_{2m-1}|&=& (\alpha-2\beta)m.
\end{array}
\end{equation}
For $\alpha<2\beta$, we have $\tau_1=\mathrm{F}(0)$, $|\Delta_1|=x_1$, and for $m\geqslant 1$
$$
\tau_{2m}= \mathrm{F}(x_{m})-x_{m},
\qquad
\tau_{2m+1}= \mathrm{F}(x_{m}')-x_{m}'
$$
and 
$$
|\Delta_{2m}|= x_{m}'-x_{m},
\hskip 30pt
|\Delta_{2m+1}|= x_{m+1}-x_{m}'
$$
giving
\begin{equation}\label{eq:ParametersII}
\begin{array}{rcll}
\tau_1&=&2\beta-\alpha\\
\tau_{2m}&=&(2m+1)(2\beta-\alpha)\\
\tau_{2m+1}&=&4m(\beta-\alpha)
&\qquad\multirow{2}{100pt}{$\alpha<2\beta, \hskip 3pt m\geqslant 1$}\\
|\Delta_1|&=&\alpha-\beta\\
|\Delta_{2m}|&=&(2m+1)(\alpha-\beta)\\
|\Delta_{2m+1}|&=&m(2\beta-\alpha).
\end{array}
\end{equation}


\bigskip
Let $\sigma$ be the permutation of $\N$ induced by $\mathrm{F}$, whereby $\sigma(j)=i$
means that the $j$th interval ends up in position $i$.

\begin{proposition} \label{prop:sigma}
{\sl The permutation $\sigma$ induced by the IET (\ref{eq:Parameters}) and (\ref{eq:ParametersII}) is given by:}
$$
\sigma(1,2,3,\ldots)=
\begin{cases}
(3,1,5,2,7,4,9,6,\ldots)&\quad \alpha>2\beta\\
(2,4,1,6,3,8,5,10\ldots)  &\quad \alpha<2\beta
\end{cases}
$$
that is, for $n=1,2,\ldots$
\begin{equation}\label{eq:Sigma}
\begin{array}{lllll}
\sigma(2)=1,&\quad&\sigma(2n+2)=2n,&\sigma(2n-1)=2n+1&\quad \alpha>2\beta\\
\sigma(1)=2,&     &\sigma(2n)=2n+2,&\sigma(2n+1)=2n-1&\quad \alpha<2\beta.
\end{array}
\end{equation}
\end{proposition}

\proof
From proposition \ref{thm:ParameterSymmetry} it suffices to consider the case $\alpha> 2\beta$. 
Note that the inverse permutation $\sigma^{-1}$ for $\alpha>2\beta$ is equal 
to the direct permutation $\sigma$ for $\alpha < 2\beta$, and vice-versa, in
agreement with proposition \ref{thm:ParameterSymmetry}.
Defining the sets of indices
\begin{equation}
L_1=\emptyset,\qquad L_i=\{\sigma^{-1}(k)\,:\, k<i\}\quad i > 1,
\end{equation}
and considering that $L_{i+1}\setminus L_i=\{\sigma^{-1}(i)\}$, we verify that
(\ref{eq:Sigma}) is equivalent to 
\begin{equation}\label{eq:L}
\begin{array}{rcl}
L_{2n+1}&=&\{1,\ldots,2(n+1)\}\setminus \{2n-1,2n+1\}\\
L_{2n}&=&\{1,\ldots,2n\}\setminus\{2n-1\}
\end{array}
\qquad n\geqslant 1.
\end{equation}

From (\ref{eq:Delta}), we have, for all $i,j$:
\begin{equation}\label{eq:IETsum}
\sigma(j)=i
\quad \iff\quad 
\mathrm{F}(\delta_{j-1})=\sum_{k\in L_i}|\Delta_k|.
\end{equation}
We shall establish the theorem via the rightmost identity, using 
formulae (\ref{eq:delta}) and (\ref{eq:Parameters}).
For $j=2$, we find
$$
\mathrm{F}(\delta_1)=\mathrm{F}(x_0')=x_0'+\tau_2=0=\sum_{k\in L_1}|\Delta_k|
$$
(the sum is empty) which establishes that $\sigma(2)=1$.

Next we let $j=2n-1$, and we shall use the identity
\begin{equation*}\label{eq:diff1}
x_{n+1}-x_{n-1}=\alpha(2n+1)-2\beta\qquad n\geqslant 1,
\end{equation*}
derived from (\ref{eq:xmyn}).
We compute:
\begin{eqnarray*}
\mathrm{F}(\delta_{2n-2})&=&\mathrm{F}(x_{n-1})=x_{n-1}+\tau_{2n-1}\\
&=&x_{n+1}-\alpha(2n+1)+2\beta+4\beta n\\
&=&x_{n+1}-(\alpha-2\beta)(2n+1)\\
&=&\delta_{2(n+1)}-|\Delta_{2n-1}|-|\Delta_{2n+1}|\\
&=&\sum_{k=1}^{2(n+1)}|\Delta_k|-|\Delta_{2n-1}|-|\Delta_{2n+1}|\\
&=&\sum_{k\in L_{2n+1}}|\Delta_k|.
\end{eqnarray*}
This shows that $\sigma(2n-1)=2n+1$, as desired.

Similarly, for $j=2n$ we need the following identity 
\begin{equation*}\label{eq:diff2}
x_{n}-x_{n}'=(2\beta-\alpha)(n+1)\qquad n\geqslant 0
\end{equation*}
derived from (\ref{eq:xmyn}) and (\ref{eq:xpm}).
Proceeding as above, we obtain:
\begin{eqnarray*}
\mathrm{F}(\delta_{2n+1})&=&\mathrm{F}(x_n')=x_n'+\tau_{2n+2}\\
&=&x_n-(\alpha-2\beta)(n+1)-(\alpha-2\beta)(2n+1)\\
&=&\delta_{2n}-(\alpha-2\beta)n\\
&=&\sum_{k=1}^{2n}|\Delta_k|-|\Delta_{2n-1}|\\
&=&\sum_{k\in L_{2n}}|\Delta_k|.
\end{eqnarray*}
This shows that $\sigma(2n+2)=2n$, and the proof is complete.
\/ $\Box$

\section{Symbolic dynamics}\label{section:SymbolicDynamics}

In accordance with the results of section \ref{section:Parameters}, 
in the rest of this paper we shall assume that $\alpha$ and $\beta$ are co-prime
and that $\alpha>2\beta$. 

We introduce several related symbolic dynamics for the interval-exchange 
transformation $\mathrm{F}$. Every $\Delta$-interval has an index $c$, given by
\begin{equation}\label{eq:c}
c(x)=n\qquad \iff\qquad x\in\Delta_n.
\end{equation}
Next we glue adjacent $\Delta$-intervals pairwise, to obtain the 
\textit{blocks} $\Xi_n$:
\begin{equation}\label{eq:Xi}
\Xi_n:=\Delta_{2n-1}\cup\Delta_{2n} \qquad n\geqslant 1.
\end{equation}
Every block has a block index $b$, given by
\begin{equation}\label{eq:b}
b(x)=n\qquad \iff\qquad x\in\Xi_n.
\end{equation}
Thus
\begin{equation}\label{eq:c2b}
b(x)=\left\lfloor \frac{c(x)+1}{2}\right\rfloor.
\end{equation}

The \textit{code} $C(x)=(c_0,c_1,c_2,\ldots)$ of a point $x\in\Z_+$ is the 
sequence of natural numbers that label the intervals visited by the orbit 
of $x$, that is, $c_t=c(\mathrm{F}^t(x))$, with $c$ given by (\ref{eq:c}).
The \textit{block code} $B(x)=(b_0,b_1,\ldots)$ is defined similarly,
using the function (\ref{eq:b}).

We shall also consider translated codes, using the notation
\begin{equation}\label{eq:ShiftedCode}
C(x)+k:=(c_0+k,c_1+k,c_2+k,\ldots).
\end{equation}
The \textit{minimum point} $\eta(x)$ is the smallest element
of the orbit through $x$, namely
\begin{equation}\label{eq:mu}
\eta(x):=\min\{\mathrm{F}^t(x)\,:\,t\in\Z \}.
\end{equation}
The \textit{transit time} $t_\eta(x)$  is defined to be the integer $t$
such that $\mathrm{F}^t(x)=\eta(x)$, if $x$ is not periodic, and the 
smallest non-negative such integer if $x$ is periodic.
In the former case, $t_\eta$ may be negative.

We introduce two auxiliary codes, namely
\begin{equation}\label{eq:CcircCstar}
\begin{array}{rcl}
C^\circ(x)&=&C(x)-c(x)\\
C^*(x)&=&C^\circ(\eta(x))=C^\circ(\mathrm{F}^{t_\eta}(x)),
\end{array}
\end{equation}
called, respectively, the \textit{translated code} and the
\textit{normalised code} of the point $x$.
Each code defines an equivalence relation on $\Z_+$, and we shall 
denote the equivalence class of $x$ for each of the three $C$-codes by 
$[x]$, $[x]^\circ$, and $[x]^*$, respectively.

For any $x$, the set $[x]$ is a segment (by which we mean
a finite set of consecutive integers), being the intersection of 
pre-images of segments $\Delta_n$ under $\mathrm{F}$. 
On each set $[x]$, the motion is rigid.
\medskip

\begin{lemma} \label{lemma:PeriodicCode}
{\sl
The code $C(x)$ is periodic if and only if the orbit through $x$
is periodic, in which case the period of the code and that of the
orbit coincide.
}
\end{lemma}

\proof
If the code is not periodic, then the orbit cannot be periodic.
Assume now that $C(x)$ is periodic with period $T$. Since $[x]$ 
is finite, the orbit through $x$ must be periodic with period 
$nT$, for some $n\geqslant 1$. 
Now, for any $k$, we have $x-\mathrm{F}^T(x)=\mathrm{F}^{kT}(x)-\mathrm{F}^{(k+1)T}(x)$,
this difference being determined solely by the periodic part 
of the code. Thus $ 0=x-\mathrm{F}^{nT}(x)=n(x-\mathrm{F}^T(x))$
and hence $n=1$ and $x=\mathrm{F}^T(x)$. \/ $\Box$
\medskip

Given two codes $C$ and $C'$, we write $C<C'$ to mean that
either $c_0<c_0'$ or there is $i\in\N$ such that $c_i<c_i'$ and
$c_k=c_k'$ for $k=0,\ldots,i-1$. The set of all codes (of any of
the above types) is therefore totally ordered.
Using the notation (\ref{eq:ShiftedCode}), we have, for any $k$,
\begin{equation}\label{eq:ShiftedCodeOrdering}
C(x)<C(x')\quad \iff\quad C(x)+k<C(x')+k.
\end{equation}

We now let
\begin{equation}\label{eq:C_t}
C_t(x)=C(\mathrm{F}^t(x))\qquad t\in\Z
\end{equation}
be the codes for all possible initial conditions along the orbit of $x$.
Then we define the  \textit{minimum code} $\underline{C}(x)$ as
\begin{equation}\label{eq:MinimumCode}
\underline{C}(x)=\min \{C_t(x)\,:\,t\in\Z\}
\end{equation}
where the minimum is computed with respect to the above ordering.
Such a minimum obviously exists.
The following result connects the minimum point to the minimum code.

\begin{proposition}\label{prop:MinimumCodeTransitTime}
{\sl
For all $x\in\Z_+$ we have
$$
\underline{C}(x)=C_{t_\eta}(x).
$$
}
\end{proposition}
To prove this result, we need a lemma.

\begin{lemma}\label{lemma:Ordering}
{\sl 
For all $x,x'\in\Z_+$, if $C(x)<C(x')$ then $x<x'$. Conversely,
if $x<x'$, then $C(x)\leqslant C(x')$.
}
\end{lemma}

\proof
Let $C(x)<C(x')$. If $c_0<c_0'$, we have finished. Otherwise, let
$i$ be as in the definition of ordering of sequences. Since $c_i<c_i'$
we have that $\Delta_{c_i}$ lies on the left of $\Delta_{c_i'}$, and 
hence $ \mathrm{F}^i(x)<\mathrm{F}^i(x')$. 
Now
$$
\mathrm{F}^i(x)=x+\sum_{k=0}^{i-1}\tau_{c_k}\,<\,
\mathrm{F}^i(x')=x'+\sum_{k=0}^{i-1}\tau_{c_k'},
$$
where the $\tau$s are the translations. By assumption, the corresponding
terms under the summation symbol are the same, and hence their
sum is the same, giving $x<x'$.

Conversely, assume that $x<x'$. If $C(x)\not=C(x')$, then
there is a smallest index $i$ for which $c_i\not= c_i'$.
If $i=0$, then $c(x)<c(x')$, and we have finished.
Otherwise, the argument used above gives that $\mathrm{F}^i(x)<\mathrm{F}^i(x')$, and since
$c(\mathrm{F}^i(x))\not=c(\mathrm{F}^i(x'))$, then $c_i<c_i'$, necessarily, whence
$C(x)<C(x')$, as desired.
\/$\Box$
\bigskip

\noindent{\sc Proof of proposition \ref{prop:MinimumCodeTransitTime}.}
Let $C_t(x)$ be as in (\ref{eq:C_t}).
We will show that $C_{t_\eta}(x)\leqslant C_t(x)$, for all $t\in\Z$.
Let $c^-$ be the smallest code element:
$$
c^-(x)=\min\{c_t(x)\,:\,t\in\Z\},
$$
and let
$$
T(x)=\{t\in\Z\,:\,c_t(x)=c^-(x)\}.
$$
Clearly, $t_\eta\in T(x)$. 
If $t\not\in T(x)$, then $C_{t_\eta}(x)<C_{t}(x)$, since 
the former code has a smaller first element. 
So we only need to show that $C_{t_\eta}(x)<C_t(x)$ for 
$t\in T(x)\setminus\{t_\eta\}$.
For this purpose it suffices to establish that $C_{t_\eta}(x)\not=C_t(x)$. Indeed, if
$C_{t_\eta}(x)$ were greater than $C_t(x)$, then lemma \ref{lemma:Ordering} 
would give $\mathrm{F}^{t_\eta}(x)>\mathrm{F}^t(x)$, contrary to the definition of minimum point.

We have two cases.

Case I: $T(x)$ is finite. Then the orbit through $x$ is not periodic.
Take any $t\in T(x)\setminus \{t_\eta\}$. Then the number of entries 
$c^-$ appearing in the codes $C_{t_\eta}(x)$ and $C_t(x)$ is different, 
and hence $C_{t_\eta}(x)\not=C_t(x)$, as desired.

Case II: $T(x)$ is infinite. Then the orbit is periodic, since $\mathrm{F}$ is invertible
and the orbit visits infinitely many times the finite set $\Delta_{c^-}$. 
Let $\ell$ be the period of the orbit (hence of the code, from lemma 
\ref{lemma:PeriodicCode}), and choose $t$ in the range $t_\eta<t<t_\eta+\ell$ .
Then the quantity $\delta=x_t-x_{t_\eta+\ell}$ is positive, because $x_{t_\eta+\ell}$ is
the minimum point and $x_t$ is not. Assume now that $C_{t_\eta}(x)=C_t(x)$.
Then $x_{t_\eta}-x_{t_\eta+\ell-t}$ is also equal to $\delta$, since it is determined by 
the same code. But this would imply that $x_{t_\eta+\ell-t}$ is smaller than
the minimum point, a contradiction.  Thus $C_{t_\eta}(x)\not=C_t(x)$, as desired.
The proof is complete.
\/ $\Box$

\section{The reduced system}\label{section:ReducedSystem}

If we order the cylinder sets $[\,\cdot\,]$ of the $C$-code according to 
the lexicographical ordering, then from lemma \ref{lemma:Ordering}, the resulting 
sequence $X_0,X_1,\ldots$, has $X_0=[0]$, and $X_{n+1}$ lying immediately 
to the right of $X_n$. The dynamics of $\mathrm{F}$ on $\mathbb{Z}_+$ induces a 
dynamics on cylinder sets $[x]\mapsto [\mathrm{F}(x)]$, which we shall
represent as dynamics on integers.
There are two problems to be dealt with. First, there are anomalies near 
the origin; these are circumvented by looking at large amplitudes.
Second, there are anomalous cylinder sets, whose size does not increase linearly with
the block order; these are dealt with by scaling.

In this section we derive the so-called \textit{reduced} interval-exchange map $\mathrm{F}'$, 
obtained from $\mathrm{F}$ by scaling coordinates in such a way as to
obtain a spatially periodic system, whose period is the block size.
The points in the phase space of the reduced system represent the
so-called \textit{regular} cylinder sets of the $C$-code.
The latter correspond to a set of full measure of orbits of $\mathrm{F}$, as
we shall see in section \ref{section:RegularPoints}.
\medskip

Using formulae (\ref{eq:Parameters}), the following asymptotic
relations for $m\to\infty$ are established at once:
\begin{equation}\label{eq:AsymptoticParameters}
\begin{array}{rclcrcl}
|\Delta_{2m-1}|&\sim& (\alpha-2\beta)m
&\qquad&
|\Delta_{2m}|&\sim& 2\beta m\\
\tau_{2m+1}&\sim& 4\beta m
&&
\tau_{2m}&\sim& 2(2\beta-\alpha) m.
\end{array}
\end{equation}

Let $\Xi_m$ be as in (\ref{eq:Xi}). Then
\begin{equation}\label{eq:SizeOfBlocks}
|\Xi_m|=|\Delta_{2m-1}|\cup|\Delta_{2m}|=\alpha m -\beta,
\qquad
|\Xi_m|\sim \alpha m.
\end{equation}

\begin{figure}[t]
\centering
\begin{tikzpicture}[scale=0.35]
\draw [blue, thick] (0,0) -- (5,2) -- (7,0) -- (12,2) -- (14,0) -- (19,2) -- (21,0) -- (26,2) -- (28,0) -- (30,0.8);
\draw [blue, thick] (0,-6.8) -- (2,-6) -- (4,-8) -- (9,-6) -- (11,-8) -- (16,-6) -- (18,-8) -- (23,-6) -- (25,-8) -- (30,-6);
\draw [thick] (2.2,0.6) rectangle (3.2,1.6);
\draw [thick] (6.0,-6.5) rectangle (7.0,-7.5);
\draw [thick] (13,1) circle (0.5);
\draw [thick] (3,-7) circle (0.5);
\draw [green, dotted, very thick] (-1,-3) -- (31,-3);
\draw [green,very thick] (0,-3.25) -- (0,-2.75);
\draw [green,very thick] (7,-3.25) -- (7,-2.75);
\draw [green,very thick] (14,-3.25) -- (14,-2.75);
\draw [green,very thick] (21,-3.25) -- (21,-2.75);
\draw [green,very thick] (28,-3.25) -- (28,-2.75);
\draw [->, ultra thick] (20,-5) -- (20,-1);
\draw [ultra thick] (20,-1) -- (20,1);
\draw [->, ultra thick] (10,-7) -- (10,-1);
\filldraw[black] (20,-3) circle (8pt);
\node  at (19.2,-2.3) {$x$};
\filldraw[black] (10,-3) circle (8pt);
\node  at (8.6,-2.3) {$\mathrm{F}(x)$};
\filldraw[black] (12,2) circle (8pt);
\filldraw[black] (26,2) circle (8pt);
\filldraw[black] (2,-6) circle (8pt);
\filldraw[black] (16,-6) circle (8pt);
\filldraw[black] (30,-6) circle (8pt);
\node  at (12,3) {$a$};
\node  at (26,3) {$a$};
\node  at (2,-5) {$a$};
\node  at (16,-5) {$a$};
\node  at (30,-5) {$a$};
\filldraw[red] (0,0) circle (8pt);
\filldraw[red] (14,0) circle (8pt);
\filldraw[red] (28,0) circle (8pt);
\filldraw[red] (4,-8) circle (8pt);
\filldraw[red] (18,-8) circle (8pt);
\node at (0,-1) {$b$};
\node at (14,-1) {$b$};
\node at (28,-1) {$b$};
\node at (4,-9) {$b$};
\node at (18,-9) {$b$};
\end{tikzpicture}
\caption{\label{fig:TranslationSurface}\rm\small
The translation surface of the reduced system, constructed from 
the infinite region lying between the two polygonal lines, by identifying pairs 
of parallel sides according to (\ref{eq:ReducedIET}). 
(Two pairs of corresponding sides are marked explicitly.)
The IET is the first-return map to the dotted line for the vertical flow, 
the ticks marking the boundary of the blocks. 
The points $a$ and $b$ are two of the four infinitely branched singular points 
on the surface. 
}
\end{figure}
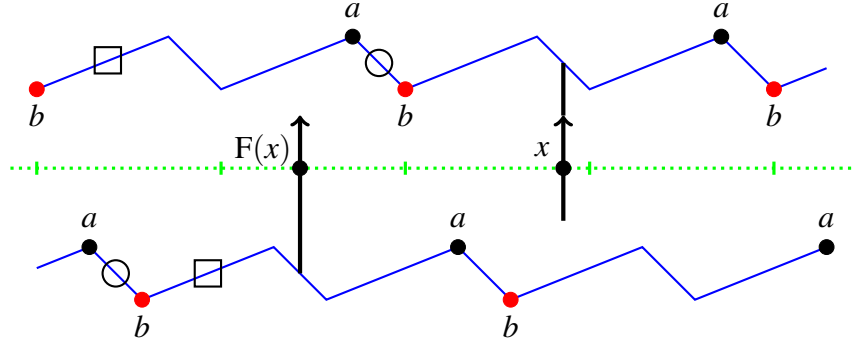

Scaling by $m$, and taking the limit $m\to\infty$, we obtain a periodic 
interval-exchange transformation, whose period is the block length $\alpha$, 
which we then extend to the whole of $\mathbb{Z}$. 
For definiteness, we shall place the left end-point of the interval 
$\Delta_1$ at the origin. (We shall make a different choice in section 
\ref{section:RegularPoints}.)
This is the reduced system.
For any $m\in\mathbb{Z}$, we have:
\begin{equation}\label{eq:ReducedIETdata}
\begin{array}{rclcrcl}
|\Delta_{2m-1}'|&=& \alpha-2\beta
&\qquad&
|\Delta_{2m}'|&=& 2\beta\\
\tau_{2m-1}'&=& 4\beta
&\qquad&
\tau_{2m}'&=& 2(2\beta-\alpha)
\end{array}
\qquad \alpha>2\beta
\end{equation}
from which we obtain
\begin{equation}\label{eq:ReducedIET}
\mathrm{F}':\mathbb{Z}\to\mathbb{Z}
\hskip 40pt
z\mapsto\begin{cases} 
z+4\beta & \mbox{if}\qad z\mod{\alpha}<\alpha-2\beta\\
z-2(\alpha-2\beta)& \mbox{otherwise.}
\end{cases}
\end{equation}
The reduction of $\mathrm{F}'$ modulo $\alpha$ is a rotation:
\begin{equation}\label{eq:Fmodalpha}
{\mathrm{F}'}(z)\equiv z+4\beta \mod{\alpha}.
\end{equation}
The translation surface of $\mathrm{F'}$ is depicted in figure \ref{fig:TranslationSurface}.

The $B$ and $C$-codes for the reduced system are defined in the obvious way.
Then we determine the domains corresponding to transitions between intervals
and blocks.
Four distinct parameter ranges need to be considered.
In each case, we display a partition of the blocks consisting of four half-open intervals.
We provide the length of each interval, and two associated transitions: 
\begin{equation}\label{eq:DifferenceCodes}
\mathrm{d} c(z)=c(\mathrm{F}'(z))-c(z)
\hskip 40pt
\mathrm{d} b(z)=b(\mathrm{F}'(z))-b(z).
\end{equation}
The former is the transition between IET domains, expressed as the change of the 
$c$-code for both odd-order (1) and even order (2) intervals; the latter is
the transition between blocks, expressed as the change of the $b$-code.

\bigskip

\goodbreak\noindent Case I: $0\leqslant 6\beta<\alpha$

\begin{equation}\label{eq:Transitions1}
\begin{array}{lllr}
\mbox{interval}&\mbox{length}&\qad\mathrm{d} c&\mathrm{d} b\\
\noalign{\vskip 10pt}
0\leqslant z< \alpha-6\beta&  \alpha-6\beta &1:0&0\\
\alpha-6\beta\leqslant z< \alpha-4\beta& 2\beta &1:+1&0\\
\alpha-4\beta\leqslant z< \alpha-2\beta& 2\beta &1:+2&+1\\
\alpha-2\beta\leqslant z< \alpha& 2\beta &2:-3&-1
\end{array}
\end{equation}
\bigskip

\goodbreak\noindent Case II:  $4\beta\leqslant\alpha<6\beta$

\begin{equation}\label{eq:Transitions2}
\begin{array}{lllr}
0\leqslant z< \alpha-4\beta& \alpha-4\beta &1:+1&0\\
\alpha-4\beta\leqslant z< \alpha-2\beta& 2\beta &1:+2&+1\\
\alpha-2\beta\leqslant z< 2(\alpha-3\beta)& \alpha-4\beta &2:-3&-1\\
2(\alpha-3\beta)\leqslant z< \alpha& 6\beta-\alpha &2: -2&-1
\end{array}
\end{equation}
\bigskip

\goodbreak\noindent Case III:  $3\beta\leqslant\alpha<4\beta$

\begin{equation}\label{eq:Transitions3}
\begin{array}{lllr}
0\leqslant z< 2(\alpha-3\beta)& 2(\alpha-3\beta) &1:+2&+1\\
\alpha-4\beta\leqslant z< \alpha-2\beta& 4\beta-\alpha &1:+3&+1\\
\alpha-2\beta\leqslant z< 2(\alpha-2\beta)& \alpha-2\beta &2:-2&-1\\
2(\alpha-2\beta)\leqslant z< \alpha& 4\beta-\alpha &2:-1&0
\end{array}
\end{equation}
\bigskip

\goodbreak\noindent Case IV:  $2\beta\leqslant\alpha<3\beta$

\begin{equation}\label{eq:Transitions4}
\begin{array}{lllr}
0\leqslant z< \alpha-2\beta& \alpha-2\beta &1:+3&+1\\
\alpha-2\beta\leqslant z< 2(\alpha-2\beta)& \alpha-2\beta &2:-2&-1\\
2(\alpha-2\beta)\leqslant z< 3(\alpha-2\beta)& \alpha-2\beta &2:-1&0\\
3(\alpha-2\beta)\leqslant z< \alpha& 2(3\beta-\alpha) &2: 0&0
\end{array}
\end{equation}
\bigskip

From the above data, we see that the $c$-code can be recovered from the $b$-code, as follows:
\begin{equation}\label{eq:b2c}
\begin{array}{rclcl}
c(z)&=&\begin{cases}2b(z)&\mbox{if}\qad \mathrm{d}b(z)=-1\\2b(z)-1&\mbox{otherwise}\end{cases}
     &\hskip 40pt&4\beta\leqslant \alpha\\
\noalign{\vskip 3pt}
c(z)&=&\begin{cases}2b(z)-1&\mbox{if}\qad \mathrm{d}b(z)=+1\\2b(z)&\mbox{otherwise}\end{cases}
     &\hskip 40pt&4\beta>\alpha.
\end{array}
\end{equation}
Thus $B$ determines $C$, while the inverse relation is established by (\ref{eq:c2b}). 

With reference to (\ref{eq:DifferenceCodes}), 
we let the difference code $D(z)$ (or $d$-code) be defined as follows:
$$
D(z)=(d_1,d_2,\ldots)\qquad d_t(z)=\mathrm{d}b(z)=b((\mathrm{F'})^{t+1}(z))-b((\mathrm{F'})^t(z)).
$$
Clearly $D(z)$ is a re-coding of $B(z)$, hence of $C(z)$.

The next result establishes the dynamics of the reduced system.

\begin{teorema}\label{thm:ReducedSystem}
{\sl Let $\overline\alpha={\alpha}/{\gcd(\alpha,2\beta)}$.
If $\overline\alpha$ is odd, then all orbits of the reduced system are periodic
with period $\overline\alpha$; in addition, all orbits have the same normalised code. 
If $\overline\alpha$ is even, then all orbits escape to $\pm\infty$. 
Specifically, if we stipulate that 0 is the left end-point of a block, then
\begin{equation}\label{eq:epsilon}
(\mathrm{F}')^{\alpha/4}(z)=z+\epsilon(z)\alpha\hskip 40pt
\epsilon(z)=\begin{cases}
+1 &\mbox{if \/ } z\equiv 0,1 \mod{4}\\
-1 &\mbox{if \/ } z\equiv 2,3 \mod{4}.
\end{cases}
\end{equation}
}
\end{teorema}

\proof 
We consider the transition domains with non-zero value of $\mathrm{d}b$. 
From tables (\ref{eq:Transitions1})-(\ref{eq:Transitions4}) we see that for any choice 
of parameters, the interval with 
$\mathrm{d}b=+1$ and that with $\mathrm{d}b=-1$ have the same length.

If $\overline\alpha$ is odd, then we distinguish two cases.
If $\alpha$ is odd, then there is a single orbit modulo $\alpha$. Because the transition
intervals have equal length, we have
\begin{equation}\label{eq:ZeroVariation}
\kappa(z)=\sum_{k=0}^{\alpha-1} \mathrm{d}b(z_t)=0,
\end{equation}
that is, modular periodicity corresponds to periodicity in $\Z$. If $\alpha$ is even,
then there are two orbits of period $\overline\alpha=\alpha/2$, and the transition
intervals have even length. Since each orbit has the same number of elements in each 
interval, equation (\ref{eq:ZeroVariation}) holds as well. Furthermore, both orbits
have the same code.

If $\overline\alpha$ is even, then $\alpha$ is divisible by 4. The number of elements of the 
two transition intervals with non-zero value of $\mathrm{d}b(z)$ is divisible by 2 but not by 4. 
Furthermore such intervals are adjacent, and their combined length is divisible by 4.
In the dynamics modulo $\alpha$ there are four orbits of period $\alpha/4$, from which it follows
that the sum $\kappa(z)$ is equal to $+1$ for two orbits and to $-1$ for the other two. 
Inspecting formulae (\ref{eq:Transitions1})--(\ref{eq:Transitions4}), we see that for all 
parameter ranges the left end-point of the $\mathrm{d}b=+1$ region is congruent modulo 4
to the left end-point of the block. Considering that the length of that region is
congruent to 2 modulo 4, if we place the origin at the left end-point of the first block, 
then it follows that $\kappa(z)=+1$ if $z\equiv0,1\mod{4}$ and $-1$ otherwise,
which is formula (\ref{eq:epsilon}).
\quad $\Box$

\section{Regular points}\label{section:RegularPoints}

We define the \textit{$\alpha$-code} of a point $x$ to be the finite sequence 
consisting of the first $\alpha$ terms in the code $C(x)$ under $\mathrm{F}$.
An $\alpha$-code of $\mathrm{F}$ is said to be \textit{regular} if it is
also the $\alpha$-code of some orbit of $\mathrm{F}'$. In this context,
we also use the terms \textrm{regular point} (a point whose $\alpha$-code is regular),
\textrm{regular cylinder set} (the cylinder set of a regular $\alpha$-code), etc.

Plainly, irregular points must exist, because the phase space of $\mathrm{F}$ is bounded 
below and that of $\mathrm{F}'$ is not. Moreover, if $\overline\alpha$ is even, then
the number of irregular points is necessarily infinite, since there is an infinite 
number of orbits with a minimum point.
The situation far from the origin is captured by the following conjecture.

\bigskip \noindent {\bf Conjecture.} 
{\sl 
If $\overline{\alpha}$ is odd, then all but finitely many points are regular.
If $\overline\alpha$ is even, then all but finitely many blocks have the same 
(positive) number of irregular points.
}
\medskip

In this section we establish the following weaker statement.
\begin{teorema} \label{thm:FullDensity}
{\sl
Let $\Gamma$ be the set of regular points of the Poincar\'e map $\mathrm{F}$.
Then $\Gamma$ has full density. Moreover, we have the block decomposition
$$
\Xi_n=\bigcup_{k=0}^{\overline\alpha-1}\Xi_{n,k}+\Lambda_n
$$
where the $\Xi_{n,k}$s are the regular cylinder sets in the $n$th block,
ordered from left to right (equivalently, by the code ordering
introduced in section \ref{section:SymbolicDynamics}), and
$$
|\Xi_{n,k}|=n\,\gcd(\alpha,2\beta)+O(1),
\hskip 40pt
|\Lambda_n|=O(1).
$$
}
\end{teorema}

\proof
We fix a sufficiently large integer $n$.
Let the \textit{stretched map} be the IET obtained from $\mathrm{F}'$ 
[see (\ref{eq:AsymptoticParameters})] by 
multiplying by $n$ all interval lengths and translations. 
In what follows, the symbols $\mathrm{F}',\Xi',\Delta',\tau'$ will refer 
to the stretched system for the current choice of $n$, with 
$\Delta'_n=[\delta'_{n-1},\delta'_n)$ [cf.~(\ref{eq:delta}) and 
(\ref{eq:Delta})].

The left and middle singularities of the block $\Xi_{n+k}$ are, respectively,
$$
\delta_{2(n+k-1)}=x_{n+k-1}
\hskip 40pt
\delta_{2(n+k)-1}=x'_{n+k-1}.
$$
We align the left end-points of the blocks $\Xi_n$ and $\Xi'_n$ by letting 
$\delta'_{2(n-1)}=\delta_{2(n-1)}=x_{n-1}$. 
Then the left and middle singularities of the block $\Xi'_{n+k}$ are, respectively,
$$
\delta'_{2(n+k-1)}=x_{n-1}+kn\alpha
\hskip 40pt
\delta'_{2(n+k)-1}=x_{n-1}+kn\alpha+n(\alpha-2\beta).
$$
The mismatch of the corresponding singularities is, respectively,
\begin{eqnarray*}
\partial\delta^{\ell}_{n,k}&=&\delta_{2(n+k-1)}-\delta'_{2(n+k-1)}=
 \frac{1}{2} k(\alpha k-\alpha-2\beta)
\\
\partial\delta^{m}_{n,k}&=&\delta_{2(n+k)-1}-\delta'_{2(n+k)-1}=
\frac{1}{2}k(\alpha k+\alpha-6\beta),
\end{eqnarray*}
is independent of $n$. 
Hence the quantity
$$
b_1=\max_{|k|\leqslant \alpha}\{|\partial\delta^{\ell}_{n,k}|,
  |\partial\delta^{m}_{n,k}|\}
$$
which represents the maximum distance between singularities over the 
largest region that can be spanned with $\alpha$ iterates, is 
independent of both $n$ and $k$.

The difference between the corresponding translations are given 
by\footnote{Consider that $\tau'_{n+k}=\tau'_n$.}
\begin{eqnarray*}
\partial\tau_{2(n+k)-1}&=&\tau_{2(n+k)-1}-\tau'_{2(n+k)-1}=4\beta k\\
\partial\tau_{2(n+k)}&=&\tau_{2(n+k)}-\tau'_{2(n+k)}=(2\beta-\alpha)(2k-1).
\end{eqnarray*}
We choose $x\in\Xi'_n$, and we let $C'(x)=(c'_0,c'_1,\ldots,c'_{\alpha-1})$ be the 
$\alpha$-code of $x$ under $\mathrm{F}'$.
The data (\ref{eq:Transitions1})--(\ref{eq:Transitions4})
show that the orbit of $x$ under $\mathrm{F}'$ will sweep at most $\alpha$ 
adjacent blocks, so that $k$ will be in the range
$[n-\alpha+1,n+\alpha-1]$.

The maximum distance between the orbit of $\mathrm{F}'$ and the orbit 
of $\mathrm{F}$ with the same code is estimated as follows:
\begin{eqnarray*}
\max_{t\leqslant \alpha} \bigl|\sum_{i=0}^{t-1}\partial \tau_{c'_i}\bigr|
&\leqslant& \max_{t\leqslant\alpha}\sum_{i=0}^{t-1}|\partial \tau_{c'_i}|
\leqslant\max_{t\leqslant\alpha} \bigl(t\max_{i\leqslant\alpha}
|\partial \tau_{c'_i}|\bigr)\\
&\leqslant&\alpha\max_{i\leqslant \alpha}|\partial \tau_{c'_i}|=:b_2.
\end{eqnarray*}
Since $\partial \tau_{c'_i}$ is independent of $n$, so is the 
constant $b_2$.

Let now $b=b_1+b_2$; then $b$ depends on $\alpha$ and $\beta$ but not on $n$.
Let us choose $n>2b$. Now, the equivalence classes of the stretched 
system coincide with the $\alpha$-cylinder sets. From the argument used in 
the proof of theorem \ref{thm:ReducedSystem} we deduce that each class has size 
$n\gcd(\alpha,2\beta)$. 
It follows that for each class $[y]'\subset\Xi'_n$ we can find a point $x$ 
which lies at distance greater than $b$ from the end-points of $[y]'$. 
Now consider the first $\alpha$ points in the orbit of $x$ 
under the maps $\mathrm{F}$ and $\mathrm{F}'$, respectively. Because of the way $b$ was defined,
no singularity of $\mathrm{F}$ or $\mathrm{F}'$ will lie between corresponding points
of the two orbits. This means that the $\alpha$-codes of the two maps are the same.
Since the same code is clearly available for the unstretched system, we have
that the point $x$ belongs to some regular cylinder set $\Xi_{n,k}\subset \Xi_n$,
and that such set has size $n\gcd(\alpha,2\beta)+O(1)$, where $O(1)<2b$.

Let $\Gamma$ be the union of all regular $\alpha$-cylinder sets.
Keeping in mind that $|\Xi_n'|-|\Xi_n|=\beta$, we have
$$
|\Gamma \cap \Xi_n|\geqslant \alpha(n-2b)-\beta 
$$
and 
$$
|\Lambda_n|=|\Xi_n\setminus \Gamma|\leqslant 2\alpha b.
$$
The density $\mathcal{D}(\Gamma)$ of $\Gamma$ is then given by
\begin{eqnarray*}
\mathcal{D}(\Gamma)=\lim_{N\to\infty}\frac{1}{N}\{x\in \Gamma\,:\,x\leqslant N\}
&=&\lim_{n\to\infty}\frac{1}{x_n}\sum_{k=1}^n\bigl|\Gamma\cap \Xi_k\bigr|\\
&\geqslant&\lim_{n\to\infty}\frac{1}{x_n}\sum_{k=2b+1}^n \bigl(\alpha(n-2b)-\beta\bigr)\\
&=&\lim_{n\to\infty}\frac{1}{x_n}\left[\frac{\alpha n(n+1)}{2}+O(n)\right]=1,
\end{eqnarray*}
where we have used the expression (\ref{eq:xmyn}) for $x_n$. 
This is the desired result. \hfill $\Box$
\bigskip

\section{Periodic orbits}\label{section:PeriodicOrbits}

In this section we prove the first statement of theorem \ref{thm:Main}: 
if $\overline{\alpha}$ is odd, then the periodic points have have full density.

Let $C$ be a regular $\alpha$-code, and let $C_0$ and $C_1$ be, respectively, the 
multi-sets of even and odd integers in $C$.
From the periodicity of the reduced orbit and (\ref{eq:ReducedIETdata}), we find:
\begin{equation}\label{eq:PeriodicityCondition}
0=\sum_{c\in C}\tau'_{c}=4\beta|C_1|+2(2\beta-\alpha)|C_0|,
\end{equation}
and since $|C_1|+|C_0|=\alpha$ 
(if $\alpha$ is even, we go through the period twice), 
we have
\begin{equation}\label{eq:C0C1}
|C_1|=\alpha-2\beta
\qquad
|C_0|=2\beta.
\end{equation}
Let us now consider an $\mathrm{F}$-orbit driven by the same code.
Using (\ref{eq:Parameters}) and (\ref{eq:C0C1}) we obtain:
\begin{eqnarray}
\sum_{c\in C}\tau_{c}&=&\sum_{c\in C_0} \tau_c+\sum_{c\in C_1} \tau_c\nonumber \\
&=&\sum_{c\in C_0}(2\beta-\alpha)(c-1)+\sum_{c\in C_1}2\beta(c+1)\nonumber \\
&=&(2\beta-\alpha)\sum_{c\in C_0}c+2\beta\sum_{c\in C_1}c+(\alpha-2\beta)|C_0|+2\beta|C_1|\nonumber \\
&=&(2\beta-\alpha)\sum_{c\in C_0}c+2\beta\sum_{c\in C_1}c+4\beta(\alpha-2\beta).\label{eq:Discrepancy}
\end{eqnarray}

An $\mathrm{F}$-orbit with $\alpha$-code $C$ will be periodic iff the rightmost expression is zero.
We begin to analyse this expression by introducing the following function:
\def\myskip{\hskip -5pt}
\begin{equation}\label{eq:S}
S:\mathbb{Z}\to\mathbb{Z}\hskip 30pt
S(x)=(\alpha-2\beta)\myskip\sum_{c\in C_0(x)}\myskip c-2\beta\myskip \sum_{c\in C_1(x)}\myskip c
\end{equation}
where $C=C(x)$ is the $\alpha$-code of the reduced system $\mathrm{F}'$,
and $C_0$ and $C_1$ are the multi-sets of even and odd elements in $C$.

\begin{lemma} \label{lemma:S}
{\sl
If $\overline{\alpha}$ is odd, then the function $S$ is constant.
}
\end{lemma}

\proof
Since $\overline{\alpha}$ is odd, every $\mathrm{F}'$-orbit is periodic with 
period $\alpha/\gcd(\alpha,2\beta)$, from theorem \ref{thm:ReducedSystem}.
Since the ordering of the elements of $C=C(x)$ is immaterial,
the value of $S=S(x)$ is the same for all points of the orbit of $x$.
Now, for any integer $k$ we have, using (\ref{eq:C0C1}):
\def\mskp{\hskip -7pt}
\begin{eqnarray*}
(\alpha-2\beta)\sum_{c\in C_0+k}\mskp c-2\beta\sum_{c\in C_1+k}\mskp c&=&S+
(\alpha-2\beta)\sum_{c\in C_0}k-2\beta\sum_{c\in C_1}k\\
&=&S+k(\alpha-2\beta)|C_0|-2\beta k|C_1|\\
&=&S+k(\alpha-2\beta)2\beta-2\beta k(\alpha-2\beta)=S.
\end{eqnarray*}
It follows that the value of $S$ is the same if we replace $C$ with
the normalised code $C^*$. Theorem \ref{thm:ReducedSystem} says
that there is only one normalised code. Hence, in the periodic case,
$S$ is constant.  \hfill $\Box$
\bigskip

We define the analogue of $S$ for the $b$-code:
\begin{equation}\label{eq:R}
 R(x)=(\alpha-2\beta)\sum_{c\in C_0(x)}b(c)-2\beta\sum_{c\in C_1(x)}b(c)
\end{equation}
where the sum is taken over the codes of the first $\alpha$ points of the orbit
with initial condition $x$.
The expressions $S$ and $R$ are related as follows:
\begin{eqnarray}
S&=& (\alpha-2\beta)\sum_{c\in C_0}2b-2\beta\sum_{c\in C_1}(2b-1)\nonumber \\
 &=& 2\bigl[(\alpha-2\beta)\sum_{c\in C_0}b-2\beta\sum_{c\in C_1}b +\beta\sum_{c\in C_1} 1\bigr]\nonumber \\
 &=& 2R+2\beta(\alpha-2\beta).\label{eq:S(R)}
\end{eqnarray}
Thus $R$ is constant, from lemma \ref{lemma:S}.

For $j=0,1$, let the set $X_j(x)$ be defined by the condition 
$x\in X_j\Leftrightarrow c(x)\in C_j$.
We define a second variant of $S$ and $R$:
\begin{equation}\label{eq:T}
 T(x)=(\alpha-2\beta)\sum_{y\in X_0(x)}y-2\beta\sum_{y\in X_1(x)}y
\end{equation}
again summing over the initial segment of an orbit with initial condition $x$. 
To express $T$ in terms of $R$, we consider quotient and remainder of the division of $y$ by $\alpha$:
$$
y=\alpha(b(y)-1)+r(y)\qquad\mbox{where}\qquad 0\leqslant r(y)<\alpha.
$$
Since $r$ is determined by the dynamics modulo $\alpha$, over an $\alpha$-segment of orbit, we have
$$
(\alpha-2\beta)\sum_{y\in X_0}1-2\beta\sum_{y\in X_1}1=(\alpha-2\beta)2\beta-2\beta(\alpha-2\beta)=0.
$$
Considering the above identity, and introducing the short-hand notation
\begin{equation}\label{eq:uw}
u=\alpha-2\beta, \qquad w=2\beta, \qquad u+w=\alpha
\end{equation}
the expression (\ref{eq:R}) with $\alpha b(y)=y+\alpha-r(y)$ gives
\begin{eqnarray*}
\alpha R(x)
&=& T(x)+\bigl[u\sum_{y\in X_0}(\alpha-r(y))-w\sum_{y\in X_1}(\alpha-r(y)\bigr]\\
&=& T(x)-u\sum_{y=u}^{\alpha-1}y+w\sum_{y=0}^{u-1}y\\
&=&T(x)-\frac{uw}{2}\alpha.
\end{eqnarray*}
So $T$ is constant as well, and, using (\ref{eq:S(R)})
\begin{equation}\label{eq:TRS}
T=\alpha R+\alpha\beta(\alpha-2\beta),
\qquad\mbox{whence}\qquad
\alpha S=2T.
\end{equation}

The following result is crucial.
\bigskip

\begin{lemma} \label{lemma:Key}
{\sl
Let $T$ be as in (\ref{eq:T}). Then, if $\overline{\alpha}$ is odd, we have
$ T=2\alpha\beta(\alpha-2\beta). $
}
\end{lemma}

\proof
We consider the uniform probability measure $\mu_0$ on the first block $\Xi_1'=(0,\ldots,\alpha-1)$, and its 
images $\mu_t$:
\begin{equation}\label{eq:Measures}
\mu_0(x)=\sum_{k=0}^{\alpha-1}\frac{1}{\alpha}\delta_{k,x}\hskip 40pt
\mu_t(x)=\mu_0((\mathrm{F}')^{-t}(x)),
\end{equation}
where $\delta_{k,x}$ is Kronecker's delta.
For $j=0,1$, let $\chi_j$ be the characteristic function of the set $\{x\in \mathbb{Z}\,:\, c(x)\mod{2}=j\}$. 
We decompose $\mu_t$ as follows
$$
\mu_t(x)=\mu_{t,0}(x)+\mu_{t,1}(x)\hskip 40pt \mu_{t,j}=\mu_t(x)\chi_j(x).
$$
For all $t$, the support of $\mu_t$ consists of a complete set of residues modulo $\alpha$.
This is seen by noting that if two distinct points in the support of $\mu_t$ were congruent 
modulo $\alpha$, then they would belong to different blocks, and ---due to spatial periodicity--- 
the same would hold for their respective initial points, which is not the case. (Alternatively, 
the dynamics modulo $\alpha$ is a translation [see (\ref{eq:Fmodalpha})], for which the 
measure $\mu_0$ is invariant.) 
Since the value of $\mu_{t,j}(x)$ depends only on the value of $x$ modulo $\alpha$, 
it follows that
\begin{equation}\label{eq:mu01}
\sum_{x\in\mathbb{Z}}\mu_{t,0}(x)=\sum_{x\in\mathbb{Z}}\mu_{0,0}(x)=\frac{2\beta}{\alpha}
\hskip 40pt
\sum_{x\in\mathbb{Z}}\mu_{t,1}(x)=\sum_{x\in\mathbb{Z}}\mu_{0,1}(x)=1-\frac{2\beta}{\alpha}.
\end{equation}
Consider the random variable $\xi(x)=x$. We begin to show that the expectation 
$\mathbb{E}_t(\xi)$ with respect to $\mu_t$ does not depend on $t$.
Using the identities above, we find (all sums are over $\mathbb{Z}$):
\begin{eqnarray*}
\mathbb{E}_{t+1}(\xi)&=&\sum_{x}x\mu_{t+1}(x)=\sum_{x}x\mu_{t}((\mathrm{F}')^{-1}(x))\nonumber\\
&=&\sum_{\mathrm{F}'(y)} \mathrm{F}'(y)\mu_{t,0}(y)+\sum_{\mathrm{F}'(y)} \mathrm{F}'(y)\mu_{t,1}(y)\nonumber\\
&=&\sum_{y} (y-2u)\mu_{t,0}(y)+\sum_{y}(y+2w)\mu_{t,1}(y)\nonumber\\
&=&\sum_{y} y\bigl(\mu_{t,0}(y)+\mu_{t,1}(y)\bigr)-2u\sum_y\mu_{t,0}(y)+2w\sum_y\mu_{t,1}(y)\nonumber\\
&=&\mathbb{E}_t(\xi)-2u\frac{w}{\alpha}+2w\frac{u}{\alpha}=\mathbb{E}_t(\xi),\label{eq:Expectation}
\end{eqnarray*}
where the change in the range of summation is justified by the invertibility of $\mathrm{F}'$.
Hence 
\begin{equation}\label{eq:E0}
\mathbb{E}_t(\xi)=\mathbb{E}_0(\xi)=\frac{\alpha-1}{2}.
\end{equation}
Now we consider the evolution of the second moment
\begin{eqnarray}
\mathbb{E}_{t+1}(\xi^2)-\mathbb{E}_t(\xi^2)
&=&\sum_x x^2\mu_{t}((\mathrm{F}')^{-1}(x))-\sum_x x^2\mu_t(x)\nonumber\\
&=&\sum_x\bigl(\mathrm{F}'(x)^2-x^2\bigr)\,\bigl(\mu_{t,0}(x)+\mu_{t,1}(x)\bigr)\nonumber\\
&=&4w\sum_{x} x\mu_{t,1}(x)-4u\sum_x x\mu_{t,0}(x)\nonumber\\
&&\qquad +4w^2\sum_x\mu_{t,1}(x)+4u^2\sum_x\mu_{t,0}(x)\nonumber\\
&=&4(Q_t+uw)\label{eq:SecondMoment}
\end{eqnarray}
where
$$
Q_t=w\mathbb{E}_{t,1}-u\mathbb{E}_{t,0}
\hskip 40pt
\mathbb{E}_{t,j}=\sum_x x\mu_{t,j}(x).
$$
Using (\ref{eq:uw}), we find
\begin{eqnarray}
Q_t&=&(\alpha-u)\mathbb{E}_{t,1}-u\mathbb{E}_{t,0}
  =\alpha\mathbb{E}_{t,1}-u\bigl(\mathbb{E}_{t,0}+\mathbb{E}_{t,1}\bigr)\nonumber\\
&=&\alpha\mathbb{E}_{t,1}-u\mathbb{E}_t(\xi)=\alpha\mathbb{E}_{t,1}-u\mathbb{E}_0.
\label{eq:QI}
\end{eqnarray}
Similarly,
\begin{equation} \label{eq:QII}
Q_t=w\mathbb{E}_{t,1}-(\alpha-w)\mathbb{E}_{t,0}=-\alpha\mathbb{E}_{t,0}+w\mathbb{E}_0.
\end{equation}

Iterating (\ref{eq:SecondMoment}) over one period of the orbits of $\mathrm{F}'$, 
and using (\ref{eq:QI}), we obtain
\begin{eqnarray*}
0&=&\mathbb{E}_{\alpha}(\xi^2)-\mathbb{E}_{0}(\xi^2)
  =\sum_{t=0}^{\alpha-1}\bigl(\mathbb{E}_{t+1}(\xi^2)-\mathbb{E}_{t}(\xi^2)\bigr)\\
&=&4\alpha \bigl(uw-u\mathbb{E}_0+\sum_{t=0}^{\alpha-1}\mathbb{E}_{t,1} \bigr),
\end{eqnarray*}
which yields 
\begin{equation*}\label{eq:E1}
\sum_{t=0}^{\alpha-1}\mathbb{E}_{t,1}=\frac{u}{2}(u-w-1).
\end{equation*}
Repeating the same procedure with (\ref{eq:QII}), we find:
\begin{equation*}
\sum_{t=0}^{\alpha-1}\mathbb{E}_{t,0}=\frac{w}{2}(3u+w-1).
\end{equation*}
Combining the last two expressions and using (\ref{eq:uw}), we obtain
\begin{equation}\label{eq:E01}
(\alpha-2\beta)\sum_{t=0}^{\alpha-1}\mathbb{E}_{t,0}-
2\beta \sum_{t=0}^{\alpha-1}\mathbb{E}_{t,1}=2\alpha\beta(\alpha-2\beta).
\end{equation}

The final step is to express $T$ in terms of the sum above. 
We let $y_t^{(x)}=(\mathrm{F}')^t(x)$, and
exploit the fact that $T$ is constant, to find:
\begin{eqnarray*}
T&=&u\sum_{t=0}^{\alpha-1}y_t^{(x)}\chi_0(y_t^{(x)})-w\sum_{t=0}^{\alpha-1}y_t^{(x)}\chi_1(y_t^{(x)})\\
&=&\frac{1}{\alpha}\sum_{x=0}^{\alpha-1}
\bigl[u\sum_{t=0}^{\alpha-1}y_t^{(x)}\chi_0(y_t^{(x)})-w\sum_{t=0}^{\alpha-1}y_t^{(x)}\chi_1(y_t^{(x)}) \bigr]\\
&=&u\sum_{t=0}^{\alpha-1} \left(\sum_{x=0}^{\alpha-1}y_t^{(x)} \frac{\chi_0(y_t^{(x)})}{\alpha}\right)
-w\sum_{t=0}^{\alpha-1} \left(\sum_{x=0}^{\alpha-1}y_t^{(x)} \frac{\chi_1(y_t^{(x)})}{\alpha}\right)\\
&=&u\sum_{t=0}^{\alpha-1} \mathbb{E}_{t,0}-w\sum_{t=0}^{\alpha-1} \mathbb{E}_{t,1}.
\end{eqnarray*}
Comparison with (\ref{eq:E01}) gives the desired result. \hfill $\Box$.
\medskip

We can now complete the proof of the first statement of theorem \ref{thm:Main}.
\bigskip

\noindent {\sc Completion of the proof of the first part of Theorem \ref{thm:Main}.} \/
From lemma \ref{lemma:Key}, and the second formula in (\ref{eq:TRS}) we obtain
$$
S=\frac{2}{\alpha}T=4\beta(\alpha-2\beta).
$$
From this, and the definition (\ref{eq:S}), it follows that the total 
translation given by equation (\ref{eq:Discrepancy}) is equal to zero.
If $\gcd(\alpha,2\beta)=2$, then the code is periodic with period 
$\overline\alpha$, and hence the sum of the first $\overline\alpha$ terms
in (\ref{eq:Discrepancy}) is equal to zero.
This means that any orbit of $\mathrm{F}$ whose $\alpha$-code is the same as
some $\alpha$-code of $\mathrm{F}'$ is periodic with period $\overline{\alpha}$.
Theorem (\ref{thm:FullDensity}) now states that the density of points for which 
this property holds is 1, which completes the proof of the first statement of 
the theorem. \hfill $\Box$

\section{Escape orbits}\label{section:EscapeOrbits}

In this section we prove the second statement of theorem \ref{thm:Main}: 
if $\overline{\alpha}$ is even (hence $\alpha$ is a multiple of 4), then the 
unbounded orbits have full density.

In this parameter range all orbits of the reduced system are unbounded (theorem 
\ref{thm:ReducedSystem}), and from (\ref{eq:epsilon}) we have that 
$(\mathrm{F}')^{\alpha/4}(z)=z+\alpha\epsilon(z)$, for all $z\in\mathbb{Z}$. 
Then theorem \ref{thm:FullDensity} implies that there is a set $\Gamma$ of 
full density, such as, if $x\in\Gamma$, then $x$ has the same $\alpha$-code
as some point $z=z(x)$, and hence $\mathrm{F}^{\alpha/4}(x)$ belongs to one of 
the blocks adjacent to the block of $x$.
Moreover, the overall translation is approximately equal to the local block 
length, and we must determine its exact value [see formula (\ref{eq:DiscrepancyII})], 
to ensure that this translation can be sustained indefinitely.

Let $C$ be a regular $\alpha$-code, with $C_0$ and $C_1$ as above.
Considering the argument used in the last part of the proof of 
theorem \ref{thm:ReducedSystem}, we have
\begin{equation}\label{eq:C0C1II}
|C_0(x)|=2\beta-2\epsilon(x)\hskip 40pt |C_1(x)|=\alpha-2\beta+2\epsilon(x)
\end{equation}
so that (\ref{eq:Discrepancy}) is replaced by
\begin{equation}\label{eq:TauSum}
\sum_{c\in C(x)}\tau_{c}=-S(x)+4\beta(\alpha-2\beta)+2\epsilon(x)(4\beta-\alpha),
\end{equation}
where $S$ is defined in (\ref{eq:S}).

The functions $S,R,T$ are no longer constant. They are related by the formulae
\begin{eqnarray}
S(x)&=&2R(x)+2\beta[\alpha-2\beta+2\epsilon(x)] \label{eq:S(R)II}\\
T(x)&=&\alpha R(x)+2\alpha^2\epsilon(x)+V(x) \label{eq:T(R)II}
\end{eqnarray}
where
\begin{eqnarray*}
V(x)&=&2u[\beta-\epsilon(x)]\,[\beta-\epsilon(x)+u+1-\cos(\pi z/2)]\\
&&\quad -w\bigl[u+2\epsilon(x)\bigr]\,\bigl[\frac{u+2\epsilon(x)}{2}-1-\cos(\pi z/2)\bigr].
\end{eqnarray*}
Finally,
\begin{equation}\label{eq:S(T)II}
\alpha S(x)=2T(x)-4\alpha^2\epsilon(x)-2V(x)+\alpha w\bigl(u+2\epsilon(x)\bigr).
\end{equation}
Using (\ref{eq:C0C1II}), and keeping in mind that, for all $x$, we have
$C(x+\alpha)=C(x)+2$ and $C(x+\cos(\pi x))=C(x)$, we find
\begin{equation}\label{eq:DeltaT}
T(x)=T(x+\cos(\pi x))+2\cos(\pi x)\epsilon(x)=
T(x+\alpha)+2\alpha^2\epsilon(x).
\end{equation}

The next task is to adapt to the escape regime the probabilistic argument used in 
the periodic case (lemma \ref{lemma:Key}).
We shall require a greater generality, and consider iterates of initial measures 
supported on shifted intervals $[z,z+\alpha)$ for some $z\in\mathbb{Z}$. To lighten 
up the notation, we shall continue to use the symbol $\mu_t$ for these measures,
highlighting the dependence on $z$ only where necessary.

We decompose $\mu_t$  into the sum of $\mu_t^+$ and $\mu_t^-$, supported, respectively,
on the residue classes $0,1$ and $2,3$ modulo 4. We use the unified notation 
$\mu^\epsilon$, where $\epsilon=\pm$ refers to sign of $\epsilon(z)$ [cf.~(\ref{eq:epsilon})], 
at any point of the support of $\mu$.
We further decompose these measures into $\mu_{t,0}^\epsilon$ and $\mu_{t,1}^\epsilon$, 
corresponding to even- and odd-order intervals.
The value of $\mu_{t,j}^\epsilon(z)$ is determined by the residues of $z$ modulo $\alpha$ and
modulo 4, and hence
\begin{eqnarray*}
\sum_{z\in\mathbb{Z}}\mu_{t,0}^\epsilon&=&\sum_{z\in\mathbb{Z}}\mu_{0,0}^\epsilon=\frac{1}{2\alpha}(w-2\epsilon)\\
\sum_{z\in\mathbb{Z}}\mu_{t,1}^\epsilon&=&\sum_{z\in\mathbb{Z}}\mu_{0,1}^\epsilon=\frac{1}{2\alpha}(u+2\epsilon).
\end{eqnarray*}
We shall use the notation
$$
\E_t^\epsilon(\xi)=\sum_{z\in\mathbb{Z}}z\mu_t^\epsilon(z)
\hskip 40pt
\E_{t,j}^\epsilon(\xi)=\sum_{z\in\mathbb{Z}}z\mu_{t,j}^\epsilon(z),\quad j\in\{0,1\}.
$$
Then we have $\E_t(\xi)=\E_t^+(\xi)+\E_t^-(\xi)$. As in (\ref{eq:Expectation}), we find:
\begin{eqnarray*}
\mathbb{E}_{t+1}^\epsilon(\xi)&=&\sum_{z}z\mu_{t+1}^\epsilon(z)\\
&=&\sum_{\mathrm{F}'(y)} \mathrm{F}'(y)\mu_{t,0}^\epsilon(y)+\sum_{\mathrm{F}'(y)} \mathrm{F}'(y)\mu_{t,1}^\epsilon(y)\\
&=&\mathbb{E}_t^\epsilon(\xi)-2u\frac{w-2\epsilon}{\alpha}+2w\frac{u+2\epsilon}{\alpha}=\mathbb{E}_t^\epsilon(\xi)+2\epsilon.
\end{eqnarray*}
The above recursion relation has solution
\begin{equation}\label{eq:E_tII}
\E_t^\epsilon(\xi)=\E_0^\epsilon(\xi)+2t=\frac{1}{4}\bigl[\alpha-1+2z-2\epsilon\cos(\pi z/2)\bigr]+2\epsilon t,
\end{equation}
and a straightforward calculation gives
\begin{equation}\label{eq:deltaSecondMoment}
\mathbb{E}_\alpha^\epsilon(\xi^2)-\mathbb{E}_0^\epsilon(\xi^2)
   =\alpha^2(8+2\epsilon)+2\alpha\epsilon(2z-1)-4\alpha\cos(\pi z/2).
\end{equation}

In place of (\ref{eq:SecondMoment}) we now have, using (\ref{eq:C0C1II})
\begin{eqnarray*}
\mathbb{E}_{t+1}^\epsilon(\xi^2)-\mathbb{E}_t^\epsilon(\xi^2)
&=&4Q_t^\epsilon +\frac{2}{\alpha}\bigl[w^2u+u^2w+2\epsilon (w^2-u^2)\bigr]\\
&=&4Q_t^\epsilon+2uw+4\epsilon (w-u)\label{eq:SecondMomentII}
\end{eqnarray*}
where
\begin{equation}\label{eq:Q_tII}
Q_t^\epsilon=w\E_{t,1}^\epsilon(\xi)-u\E_{t,0}^\epsilon(\xi)
=\alpha\E_{t,1}^\epsilon(\xi)-u\E_t^\epsilon(\xi)
=-\alpha\E_{t,0}^\epsilon(\xi)+w\E_t^\epsilon(\xi).
\end{equation}

We now iterate this relation, to evaluate the telescopic sum
$\mathbb{E}_\alpha^\epsilon(\xi^2)-\mathbb{E}_0^\epsilon(\xi^2)=
\sum_{t=0}^{\alpha-1} \bigl[\mathbb{E}_{t+1}^\epsilon(\xi^2)-\mathbb{E}_t^\epsilon(\xi^2)\bigr].$
A lengthy calculation using formulae (\ref{eq:E_tII})--(\ref{eq:Q_tII}) and the 
procedure employed in the previous section gives
\begin{eqnarray}\label{eq:E01twin}
2\left(u\sum_{t=0}^{\alpha-1}\mathbb{E}_{t,0}^\epsilon-
w\sum_{t=0}^{\alpha-1}\mathbb{E}_{t,1}^\epsilon\right)
&=&
\alpha\bigl[2\alpha(\beta-2)-4\beta^2+\epsilon(1-3\alpha+8\beta-2z)
\nonumber\\
&&\qquad +2(\cos(\pi z/2)\bigr].
\end{eqnarray}

The final step is to express $T$ in terms of the above expression. Since the
functions $S,R,T$ are no longer constant, we shall need the following

\begin{lemma}\label{lemma:DeltaT4}
{\sl
If $\overline \alpha$ is even, then, for any $x$ we have \/
$
T(z+4)=T(z)-8\alpha\epsilon(z).
$
}
\end{lemma}
\proof If $\overline\alpha$ is even, then equation (\ref{eq:DeltaT}) gives 
$T(z+\alpha)-T(z)=-2\alpha^2\epsilon(z)$, so 
it suffices to show that the value of $T(z+4)-T(z)$ depends only on $\epsilon(z)$.
Introducing the notation $y_t^{(a)}=(\mathrm{F}')^t(a)$, a short calculation gives
$$
T(y_1^{(z)})-T(z)=4\alpha\epsilon(z)\bigl[u\chi_0(z)-w\chi_1(z)\bigr]=
  4\alpha\epsilon(z)\bigl[\alpha\chi_0(z)-2\beta].
$$
Let $\tau$ be the smallest positive integer $t$ such that $y_t^{(z)}\equiv z+4\mod{\alpha}$,
and let $\kappa(z)$ be defined by the equation $y_\tau^{(z)}=z+4+\alpha\kappa(z)$.
We find that $\tau\equiv \beta^{-1}\mod{\alpha/4}$, independent from $z$.
Considering that $\epsilon$ is constant along orbits, we iterate the above
relation to obtain
\begin{eqnarray*}
T(z+4)-T(z)&=&T(y_\tau^{(z)})-T(z)+2\alpha^2\epsilon(z)\kappa(z)\\
&=&2\alpha\epsilon(z)\Bigl\{\alpha\bigl[\kappa(z)+2\sum_{t=0}^{\tau-1}\chi_0(y_t^{(z)})\bigr]-4\beta\tau\Bigr\}.
\end{eqnarray*}

We must show that the expression $\kappa(z)+2\sum_{t=0}^{\tau-1}\chi_0(y_t^{(z)})$ is constant.
With references to formulae (\ref{eq:Transitions1})--(\ref{eq:Transitions4}), let $\chi^+$ and $\chi^-$ 
be the characteristic functions of the intervals defined by $\mathrm{d}b=+1$ and $\mathrm{d}b=-1$, 
respectively, and let $\chi=\chi^++\chi^-$. Then
\begin{equation}\label{eq:kappa}
\kappa(z)=\sum_{t=0}^{\tau-1}\bigl[\chi^+(y_t^{(z)})-\chi^-(y_t^{(z)})\bigr] - \zeta(z)
\qquad\zeta(z)=\delta_{b(z+4),b(z)+1},
\end{equation}
where $\delta$ is Kronecker's delta. 
Let $\alpha'=\alpha/4$; we have two cases.

Case I: $\alpha'>\beta$.
In this case  we have $\chi^-=\chi_0$, and $\chi$ is the characteristic
function of the union of intervals $[\alpha-4\beta,\alpha)+\alpha\mathbb{Z}$. 
From (\ref{eq:kappa}) we obtain 
\begin{equation}\label{eq:Constant}
\kappa(z)+2\sum_{t=0}^{\tau-1}\chi_0(y_t^{(z)})
   =\sum_{t=0}^{\tau-1}\chi(y_t^{(z)})-\zeta(z).
\end{equation}
Thus the value of the left-hand side is equal to the number of points which fall in the 
interval where $\mathrm{d}b\not=0$, decreased by one unit if $z$ and $z+4$ lie in
different blocks. We have to show that such a number is constant, with the stated exception.
By conjugating the orbit through $z_0$ for the map $X\mapsto X+4\beta\mod{\alpha}$ to 
the orbit through $z=\lfloor z_0/4\rfloor$ for the map $X\mapsto X+\beta\mod{\alpha'}$,
we reduce this problem to showing that the number of elements of set
$$
A(z)=\{z+t\beta\mod{\alpha'}\,:\,t=0,\ldots,\tau-1\}\hskip 40pt 0\leqslant z\leqslant \alpha'-1
$$ 
which lie in the interval $I_1=[\alpha'-\beta,\alpha')$, is equal to some integer $n_0$ 
for all $z\not= \alpha'-1$, and to $n_0+1$ for $z=\alpha'-1$. We introduce the
symbolic dynamics of rotation by $\beta$ on the circle $[0,\alpha')$, obtained by
assigning the symbol $0$ to the interval $I_0=[0,\alpha'-\beta)$ and the symbol $1$
to the interval $I_1$ defined above. The binary words of length $\tau$ obtained by 
varying $z$, are the same as the Sturmian words of any irrational number sufficiently 
close to $\beta$. A Sturmian language is balanced \cite[theorem 6.1.8]{PytheasFogg}, 
meaning that the number of $1$s appearing in these words assumes precisely two 
consecutive values, say, $n_0$ and $n_0+1$.

Now let 
$$
A_1(z)=A(z)\cap I_1
\hskip 40pt
N(z)=\# A_1(z).
$$
The set $A(z-1)$ is obtained from $A(z)$ by shifting all points of the latter 
to the left by one unit. The set $A(\alpha'-1)$ contains both end-points of $I_1$.
By construction, $\alpha'-1\not\in A(\alpha'-2)$, and hence, if we let 
$n_0=N(\alpha'-2)$, we have $N(\alpha'-1)=n_0+1$.
Choose $z$ such that $N(z)=n_0$. 
The only way to have $N(z-1)=n_0+1$, is that, under such a left shift,
the set $A_1(z)$ gains one point on the right, and loses no point on the left.
Then $0$ must be in $A(z)$. If $z\not=0$, then the pre-image $\alpha'-\beta$ of 
$0$ also belongs to $A(z)$, and hence one point is lost in the shift.
Thus $z=0$, namely $z-1=\alpha'-1$, as required. We have shown that there is $n_0$
such that
$$
N(z)=\begin{cases} n_0 & \mbox{if}\qad z\not=\alpha'-1\\n_0+1 & \mbox{if}\qad z=\alpha'-1.\end{cases}
$$
This means that the left-hand side of (\ref{eq:Constant}) is constant,
and hence $T(z+4)-T(z)$ depends only on $\epsilon(z)$, as required.
\smallskip

Case II: $\alpha'<\beta$. 
Then $\chi^+=\chi_1$ and $\chi$ is the characteristic function
of $[0,2(\alpha-2\beta))+\alpha\mathbb{Z}$.
The analysis is the same as that given above, with the 
opposite sign in the expression $\zeta(z)$ in (\ref{eq:kappa}). 
We shall not repeat it, for the sake of brevity.
$\Box$

\begin{lemma} \label{lemma:KeyII}
{\sl
Let $T$ and $\epsilon$ be as above. Then if $\overline{\alpha}$ is odd and
$\epsilon(z)=1$, we have
$$
T(z)=-2\alpha z+\alpha\bigl[2\beta(\alpha-2\beta)+2(4\beta-3\alpha)].
$$
}
\end{lemma}

\proof
The condition $\epsilon(z)=1$ characterises the points which escape to $+\infty$. 
Equations (\ref{eq:DeltaT}) and lemma \ref{lemma:DeltaT4} give
$$
T^+(z+\gamma)=T^+(z)-2\alpha\gamma
\hskip 40pt
\gamma(z)=2-\cos(\pi z).
$$

Using the above and lemma \ref{lemma:DeltaT4}, we obtain
\begin{eqnarray*}
T^+(z)&=&\frac{2}{\alpha}\sum_{k=0}^{\alpha/4-1}[T^+(z)+T^+(z)]\\
&=&\frac{2}{\alpha}\sum_{k=0}^{\alpha/4-1}\bigl[T^+(z+4k)+T^+(z+4k+\gamma)\bigr] +\alpha(\alpha+\gamma-4).
\end{eqnarray*}

Using lemma \ref{lemma:DeltaT4}, we compute
\begin{eqnarray*}
T^+(z)&=&u\sum_{t=0}^{\alpha-1}y_t^{(z)}\chi_0(y_t^{(z)})-w\sum_{t=0}^{\alpha-1}y_t^{(z)}\chi_1(y_t^{(z)})\\
&=&\frac{2}{\alpha}\sum_{k=0}^{\alpha/4-1}
\left[u\sum_{t=0}^{\alpha-1}y_t^{(z+4k)}\chi_0(y_t^{(z+4k)})-w\sum_{t=0}^{\alpha-1}y_t^{(z+4k)}\chi_1(y_t^{(z+4k)})\right.\\
&& \quad \left. + u\sum_{t=0}^{\alpha-1}y_t^{(z+4k+\gamma)}\chi_0(y_t^{(z+4k+\gamma)})-w\sum_{t=0}^{\alpha-1}y_t^{(z+4k+\gamma)}\chi_1(y_t^{(z+4k+\gamma)}) \right]\\
&&\quad +\alpha(\alpha+\gamma-4)\\
&=&2u\sum_{t=0}^{\alpha-1} \left(\sum_{k=0}^{\alpha/4-1}
 y_t^{(z+4k)}\frac{1}{\alpha}\chi_0(y_t^{(z+4k)})+ y_t^{(z+4k+\gamma)}\frac{1}{\alpha}\chi_0(y_t^{(z+4k+\gamma)})\right)\\
&&\quad 
-2w\sum_{t=0}^{\alpha-1} \left(\sum_{k=0}^{\alpha/4-1}y_t^{(z+4k)}\frac{1}{\alpha}\chi_1(y_t^{(z+4k)})
+ y_t^{(z+4k+\gamma)} \frac{1}{\alpha}\chi_1(y_t^{(z+4k+\gamma)})\right)\\
&&\quad +\alpha(\alpha+\gamma-4)\\
&=&2\left(u\sum_{t=0}^{\alpha-1} \mathbb{E}_{t,0}^+-w\sum_{t=0}^{\alpha-1} \mathbb{E}_{t,1}^+\right)
+\alpha(\alpha+\gamma-4).
\end{eqnarray*}

The above expressions, together with (\ref{eq:E01twin}), gives an explicit
formula for $T^+(z)$:
$$
T^+(z)=-2\alpha z+\alpha\bigl[2\beta(\alpha-2\beta)+2(4\beta-3\alpha)].
$$
The proof is complete.
$\Box$.
\bigskip

We can finally complete the proof of the second part of theorem \ref{thm:Main}.
\medskip

\noindent {\sc Completion of the proof of Theorem \ref{thm:Main}.} \/
Assume that $\overline\alpha$ is even, and let $\Gamma$ be the full density set 
specified in theorem \ref{thm:FullDensity}. 
Let $x\in \Gamma$ be given, and let us assume that the orbit of $x$ drifts to the right:
$b(\mathrm{F}^\alpha(x))=b(x)+4$. Then there are precisely two consecutive integers 
$z^*=z^*(x)$, and $z^*+1$ with the property that $z^*\equiv 0\mod{4}$ and the $\alpha$-code 
of $x$ under $\mathrm{F}$ and that of $z^*$ or $z^*+1$ under $\mathrm{F}'$ are the same. 
Lemma \ref{lemma:KeyII} and equation (\ref{eq:S(T)II}) yield
$$
S(z)=-4z+4\beta(\alpha-2\beta)+12(\beta-\alpha)-4\cos(\pi z/2).
$$
Substituting this expression in (\ref{eq:TauSum}), we finally arrive at the 
following formula for the total translation under $\alpha$ iterations of 
the Poincar\'e map $\mathrm{F}$:
\begin{equation}\label{eq:DiscrepancyII}
\mathrm{F}^\alpha(x)-x=4z+10\alpha-4\beta+4\cos(\pi z/2).
\end{equation}
One verifies that both values $z=z^*$ and $z=z^*+1$ produce the same value of the 
right-hand side of (\ref{eq:DiscrepancyII}).

Let $x$ belong to the $n$th block $\Xi_n$, that is, $x_{n-1}\leqslant x < x_n$,
with $x_n$ given by (\ref{eq:xmyn}). 
According to theorem \ref{thm:FullDensity}, the set $\Xi_n\cap \Gamma$ is partitioned 
into $\alpha/2$ regular cylinder sets $\Xi_{n,k}$, $k=0,\ldots,\alpha/2-1$, of 
$2n+O(1)$ points each, corresponding to as many distinct regular $\alpha$-codes,
plus a residual set $\Lambda_n$ of size $O(1)$,
corresponding to irregular codes. Since the point $x$ is regular and its
orbit drifts to the right, there is a unique even integer $j=j(x)$ such that $x\in\Xi_{n,j}$.
Then the points $x$ and $z^*=(n-1)\alpha+2j$ have the same $\alpha$-code for the 
maps $\mathrm{F}$ and $\mathrm{F}'$, respectively, and $z^*\equiv 0 \mod{4}$.
Substituting $z=z^*$ (or $z=z^*+1$) in (\ref{eq:DiscrepancyII}), we obtain
$$
\mathrm{F}^\alpha(x)-x=4n\alpha+6\alpha-4\beta+4(1+2j).
$$
We now compute the total translation $\Delta x$ required to move a 
point $x\in\Xi_{n,j}$ to the corresponding position within $\Xi_{n+4,j}$,
four blocks to the right.
Considering the expression (\ref{eq:SizeOfBlocks}) for the block size,
and the fact that $|\Xi_{n+k,j}|=|\Xi_{n,j}|+k+O(1)$, we obtain, for $i=4$:
$$
\Delta x=\sum_{i=0}^3\bigl[(n+i)\alpha-\beta\bigr]+4+8j=\mathrm{F}^\alpha(x)-x.
$$
This identity shows that the total translation generated by a regular $\alpha$-code 
sends a point $x\in\Xi_{n,j}$ with $\epsilon(x)=1$ into a point of $\Xi_{n+4,j}$, 
with the possible exception of $O(1)$ points at the boundary of $\Xi_{n,j}$.
Hence these translations can be sustained indefinitely.
This set of points has density 1/2, and their orbits escape to infinity. 
The result now follows from the fact that $\mathrm{F}$ is invertible, which accounts for the
escape of a complementary set of density 1/2.
$\Box$
\bigskip

For completeness, we determine $z^*(x)$ explicitly, for $x\in \Gamma$ with $\epsilon(x)=1$. 
From section \ref{section:FirstReturnMap} we find that the block $n(x)$ of $x$ is given by
$$
n(x)= \left\lfloor \frac{2\beta-\alpha + \sqrt{(\alpha-2\beta)^2+8\alpha x}}{2\alpha}\right\rfloor+1.
$$
Theorem \ref{thm:FullDensity} states that there are $2n+O(1)$ points in any regular
cylinder set of $\Xi_n$. 
Keeping in mind that the left end-point of the $n$th block is $x_n$ [see equation (\ref{eq:xmyn})]
and that the length of the $n$th block is $n\alpha-\beta$, we find
$$
z^*(x)=[n(x)-1]\alpha+\left\lfloor\alpha\,\frac{x-x_{n(x)}}{n(x)\alpha-\beta}\right\rfloor.
$$
This gives $z^*(x)=\alpha n(x)+O(1)$, and hence
$$
\mathrm{F}^\alpha(x)-x=4\alpha n(x)+O(1).
$$


\end{document}